\documentclass[12pt]{amsart}
\usepackage{amssymb}
\usepackage{amsfonts}
\usepackage{latexsym}
\usepackage{amscd}

\numberwithin{equation}{section}

\newfont{\gd}{eufm10 scaled \magstep1}
\newfont{\gs}{eufm7 scaled \magstep1}
\newfont{\gss}{eufm5 scaled \magstep1}
\newcommand{\Gbd}[1]{\mbox{\gd #1}}
\newcommand{\Gbs}[1]{\mbox{\gs #1}}
\newcommand{\Gbss}[1]{\mbox{\gss #1}}
\newcommand{\got}[1]{\mathchoice{\Gbd #1}{\Gbd #1}{\Gbs #1}{\Gbss #1}}

\newcommand{\gn}{\got n}
\newcommand{\gm}{\got m}

\newcounter{cs}
\stepcounter{cs}
\newcounter{ds}
\stepcounter{ds}

\newcommand{\casos}{\begin{itemize}}
\newcommand{\fcasos}{\end{itemize}\setcounter{cs}{1}}

\newcommand{\subcasos}{\begin{itemize}}
\newcommand{\fsubcasos}{\end{itemize}\setcounter{ds}{1}}

\newcommand{\matriu}[1]{\left(\begin{array}{#1}}
\newcommand{\fmatriu}{\end{array}\right)}

\newcommand{\ld}{{\rm lim}_{_{\kern-14pt\longrightarrow \kern3pt}}}

\newcommand{\Z}{{\mathbb{Z}}}
\newcommand{\N}{{\mathbb{N}}}
\newcommand{\Q}{{\mathbb{Q}}}
\newcommand{\R}{{\mathbb{R}}}
\newcommand{\ol}{\overline}
\newcommand{\wt}{\widetilde}

\newtheorem{lem}{Lemma}[section]
\newtheorem{corol}[lem]{Corollary}
\newtheorem{theor}[lem]{Theorem}
\newtheorem{prop}[lem]{Proposition}
\newtheorem{rema}[lem]{Remark}
\newtheorem{defi}[lem]{Definition}

\newtheorem{defidis}[lem]{Definition and Discussion}
\newtheorem{p}[lem]{}
\newtheorem{qtn}{Question}
\newtheorem{probl}[lem]{Problem}
\newtheorem*{nota}{Notation}

\begin{document}
\title[Groups with wild intervals]{Simple Riesz groups of rank one having wild intervals}
\author{Francisco Ortus, Enric Pardo}
\address{Departamento de Matem{\'a}ticas, Universidad de C{\'a}diz,Apartado 40, 11510 Puerto Real (C{\'a}diz), Spain}
\email{francisco.ortus@uca.es}\email{enrique.pardo@uca.es}
\author{Francesc Perera}
\address{Departament de Matem\`atiques, Universitat Aut\`onoma de Barcelona, 08193 Bellaterra (Barcelona), Spain.}
\email{perera@mat.uab.es}
\thanks{The first author is supported by a FPU fellow of Junta de Andaluc{\'\i}a.
The second author is partially supported by the DGI and European
Regional Development Fund through Project BFM2001--3141.
First and second authors are partially supported by PAI III grant
FQM-298 of Junta de Andaluc\'{\i}a. The third author is partially
supported by the Nuffield Foundation NUF-NAL 02, the DGI and European
Regional Development Fund through Project BFM2002-01390. The
second and third authors are also partially supported by the
Comissionat per Universitats i Recerca de la Generalitat de
Catalunya.}




\dedicatory{Gert Kjaergaard Pedersen, in memoriam} 

\begin{abstract}
We prove that every partially ordered simple group of rank one
which is not Riesz embeds into a simple Riesz group of rank one if
and only if it is not isomorphic to the additive group of the
rationals. Using this result, we construct examples of simple
Riesz groups of rank one $G$, containing unbounded intervals
$(D_n)_{n\geq 1}$ and $D$, that satisfy: (a) For each $n\geq 1$,
$tD_n\ne G^+$ for every $t< q_n$, but $q_nD_n=G^+$ (where $(q_n)$ is a sequence of relatively prime integers); (b) For every
$n\geq 1$, $nD\ne G^+$. We sketch some potential applications of
these results in the context of K-Theory.
\end{abstract}

\maketitle

\section*{Introduction}

One of the subjects of interest in the theory of partially ordered
abelian groups is the analysis of intervals, that is, non-empty,
upward directed and order-hereditary subsets. These have been used
in instances of quite different flavour. For example,
in~\cite{G-H} and~\cite{G}, they proved to be essential in
studying extensions of dimension groups. In the
papers~\cite{K0Good} (see
also~\cite{elhan}),~\cite{Perera},~\cite{Ape} and~\cite{op}, their
usage was directed towards an understanding of the non-stable
K-Theory of multipliers of $C^*$-algebras with real rank zero and
von Neumann regular rings, basically by describing the monoid of
equivalence classes of projections. Other applications can be
found in~\cite{wehrung}, where the Riesz refinement property in
monoids of intervals is studied in detail; in~\cite{WhrEq}, where
a complete description of the universal theory of Tarski's
equidecomposability types semigroups is given, and also
in~\cite{WhrCom}, as an instrument to obtain some extensions of
Edwards' Separation Theorem (see, e.g.~\cite[Theorem
11.13]{POAG}).

Since, as just mentioned, these monoids appear useful in the
context of K-Theory of operator algebras, there is a strong need
for constructing explicit examples of such monoids that help in
providing evidence towards the study of certain conditions in
multiplier algebras. In this paper we present such examples, in
the form of countable Riesz groups of rank one whose monoids of
intervals enjoy certain relevant properties, thus adding new
examples to the knowledge of Riesz groups. Our motivation for the
search of these examples can also be traced back to the following
question, asked by Goodearl in~\cite[Open Problem 30]{POAG}:
\emph{Can every partially ordered simple abelian group be embedded
in a simple Riesz group ?} This was proved to be the case by
Wehrung~(\cite{Whr3}) via a cofinal embedding. However, the fact
that part of this construction depends on Model theoretical
arguments prompts the need of finding more concrete realisations
of these type of embeddings. More concretely, the embedding result
just mentioned was used in~\cite[Example 3.14]{Whr3} to obtain an
example of a torsion-free simple Riesz group $G$ containing an
interval $D\ne G^+$ such that $2D=G^+$. Wehrung then asked whether
\emph{such example can be realized as a torsion-free Riesz group
of rank one (i.e. with positive cone isomorphic to an additive
submonoid of $\Q$)}, see~\cite[Problem 3.15]{Whr3}. This question
was answered by the second author in~\cite{P2} by constructing a
large family of simple groups that can be embedded into simple Riesz groups of rank
one~(\cite[Theorem 2.11]{P2}).

We first extend this result, by showing that \emph{an ordered
simple group of rank one $(G,G^+)$ --  which is not Riesz -- can
be embedded into a simple Riesz group of rank one if and only if
$G\not\cong\Q$}. This is done in Section~1, and used subsequently
to provide wide generalizations of Wehrung's example. The main
tool in~\cite[Example 3.14]{Whr3} is the construction of a certain
submonoid of $\Q^+$ using the submonoid of $\Z^+$ generated by $2$
and $7$. We extend this construction in Section~2 to combinations
of submonoids of $\Z^+$ generated by coprime integers $p$ and $q$,
but with considerable more extra care needed. This provides us
with an example of a simple ordered group $(G, G^+)$ that
contains a proper interval $D$, a multiple of which equals the
positive cone $G^+$. However, this is not a Riesz group. An inductive procedure, based on taking
direct limits of this type of construction, leads in Section~3 to
a first example of a Riesz group for which a whole sequence of
(proper and unbounded) intervals $(D_n)$ can be constructed; every
such interval has the property that $tD_n\neq G^+$ for every
$t\leq q_n$, and $q_nD_n=G^+$. Here, $(q_n)$ is an increasing
sequence of relatively prime integers. The inductive step is based
on a suitable amalgamation of groups (of the type considered in
Section~2, which is why we can refer to them as
basic building blocks) in a commutative diagram. A further
modification of this example, after using the refinement property
on the monoid of countably generated intervals, allows us to
achieve that the sequence of intervals is moreover decreasing.

In Section~4 we state some arithmetic results on simple
components (see below), which are used in Section~5 in order to construct an
example of a simple Riesz group $G$ containing an unbounded
interval $D\subset G^+$ such that $nD\ne G^+$ for every $n\geq 1$.
The constructions carried out in Sections~3 and~5 are combined in
Section~6 to obtain an example of a simple Riesz group that
has the properties exhibited in Sections~3 and~5.

An object which is central in this paper is that of a \emph{simple
component}, as our examples are built essentially via direct limit
constructions of simple components of various kinds. In short, a
simple component is nothing else but the group $\Z$, together with
a partial ordering that makes it simple. This, for example,
includes $(\Z,\langle k,l\rangle)$, where $\langle k,l\rangle$ is
the submonoid of $\Z$ generated by two relatively prime integers.
Simple components have been studied in different contexts, notably
with relation to K-theoretical aspects of $C^*$-algebras (see,
e.g.~\cite{Vill1},~\cite{RorVill},~\cite{EllVill}, where it is
shown that there exist simple $C^*$-algebras with stable rank one
whose $K_0$ groups are simple components). Consider also the
following question:
\begin{qtn}
\label{TomsConjecture}{\rm (\cite{Toms}) Let $N\in\N$. For
every $1\leq i\leq N$ take $q_i$ and $m_i$ in $\N$ to be
relatively prime, where $q_i$ is prime. Take moreover a positive
integer $L$ that is coprime with each $q_i$ and $m_i$. Consider
the following subsemigroup of the positive integers:
\[
S=\frac{1}{L}\left( \bigcap _{i=1}^{N}\langle q_i, m_i\rangle \right)\cap \Z\,.
\]
Can every positive cone of a simple component be expressed as $S$
for suitable choices of $N$, $(q_i)$, $(m_i)$ and $L$ ?}
\end{qtn}

The construction technique developed by Toms in~\cite{Toms}
produces, for every such $S$ as above, a simple $C^*$-algebra with
stable rank one whose $K_0$ group is isomorphic to $\Z$ with
positive cone $S$. The real rank of these examples is not zero, because otherwise they would be weakly divisible in the sense of~\cite[Section~5]{pr} (see also~\cite{ror2},~\cite{ror3}). Hence, these results suggest the
problem of constructing simple $C^*$-algebras $A$ with real rank
zero and stable rank one such that $(K_0(A),K_0(A)^+)$ is
isomorphic to one of the groups we construct in this paper (as
well as those constructed in~\cite{P1} and~\cite{P2}), by lifting
connecting maps in the direct limit expression of these groups (as
limits of simple components and order-embeddings), to
$C^*$-algebra maps between $C^*$-algebras of the type constructed
in~\cite{Toms}. Other relevant aspects of this discussion are
outlined in Section~7.

Throughout the paper we will refer to~\cite{POAG} for notations
and definitions on partially ordered abelian groups. We recall
here some basic notions that we shall use frequently. A
\emph{cone} of an abelian group $G$ is an additive submonoid $P$
of $G$ containing zero. We say that the cone $P$ is \emph{strict}
if $P\cap (-P)=\{ 0\}$. A \emph{partially ordered abelian group}
is an abelian group $G$ endowed with a strict cone $G^+$, called
the \emph{positive cone} of $G$. The usual notation for a
partially ordered group is $(G,G^+)$, and the elements of $G^+$
are referred to as the \emph{positive elements} of $G$. The order
induced by $G^+$ is denoted in this paper by $\leq_G$. We say that
$(G,G^+)$ is \emph{directed} provided that any element can be
written as a difference of two positive elements. An element $u$
in $G$ is said to be an \emph{order-unit} provided that
$0\ne u\in G^+$ and for each element $x$ in $G$, there exists $n$
in $\N$ such that $-nu\leq_G x\leq_G nu$ (note that $G$ will then be
directed). A partially ordered abelian group is said to be
\emph{simple} when it is non-zero and every non-zero positive
element is an order-unit. A partially ordered abelian group
$(G,G^+)$ satisfies the \emph{Riesz decomposition property} (or is
a \emph{Riesz group}, for short) if the following condition is met
in $G^+$: whenever $x\leq_G y_1+y_2$ in $G^+$, there exist $x_1$ and
$x_2$ in $G^+$ such that $x=x_1+x_2$ and $x_j\leq_G y_j$ for all
$j$. It is well-known that this is equivalent to the Riesz
refinement and interpolation properties (see,
e.g~\cite[Proposition 2.1]{POAG}).

If $(G,G^+)$ and $(H,H^+)$ are partially ordered abelian groups, a
\emph{positive morphism} is a group homomorphism $f\colon G\to H$
such that $f(G^+)\subseteq H^+$. A positive morphism $f\colon G\to
H$ is an \emph{order-embedding} if $f$ is one-to-one and $x\in
G^+$ whenever $f(x)\in H^+$ (in other words, $f(G^+)=f(G)\cap
H^+$).

\section{Embedding results}

In this section, we will establish some results about embedding
simple groups into simple Riesz groups, that improve those
appearing in~\cite{P2}. The first one was shown by the second
author in a rather complicated way [unpublished]. The proof we
present here was pointed out by G. Bergman.

We start by recalling some basic facts related to generalized
integers (see, e.g.~\cite{P1}). Let $\mathbb{P}$ be the set of the
natural prime numbers. A {\it generalized integer} $\gn$ is a map
$$\gn\colon\mathbb{P} \to \{0,1,2,\dots,\infty\}\,.$$ Usually we
write
\begin{eqnarray}
\label{general}
\gn =\prod _{p\in \mathbb{P}}{p^{\gn (p)}}.
\end{eqnarray}
When $\gn$ is finite (i.e. it never takes the value $\infty$ and
it is zero except at finitely many primes), we identify $\gn$ with
the integer appearing on the right hand side of~(\ref{general}).
Multiplication extends naturally to generalized integers, namely,
the product $\gn\cdot\gm$ of $\gn$ and $\gm$ is defined as
$(\gn\cdot\gm)(p)=\gn(p)+\gm(p)$ for every $p$ in $\mathbb{P}$.
Thus we say that $\gn$ divides $\gm$, in symbols $\gn\mid \gm$, if
there is $\gn'$ such that $\gm=\gn\cdot\gn'$. We say that $\gn$
and $\gm$ are coprime if for every $p$ in $\mathbb{P}$ we have
$0\in \{ \gn (p), \gm (p)\}$.

Given a generalized integer $\gn$, we associate to it an additive
subgroup of $\mathbb{Q}$ containing $1$ by setting
$\mathbb{Z}_{\gn }= \{\frac{a}{b}\in\Q\mid a\in \mathbb{Z}\mbox{
and } b\mid {\gn }\}$. Conversely, one can associate a generalized
integer to any subgroup of $\Q$ that contains $1$, and these
assignements are mutually inverse (see~\cite[Lemma 2.3]{P1}).

Given a sequence $A=(a_n)_{n\geq 1}$, we define $\gn
(A)=\prod_{n\geq 1}{a_n}$, and we say that the sequence $A$ is the
sequence associated to $\gn (A)$. A sequence $A=(a_n)$ is
associated to a generalized integer $\gn$ when $\gn=\gn(A)$. We
can always associate a sequence to a generalized integer, as shown
in~\cite[Lemma 3.10]{P1}.

One notion that will become relevant in this paper is that of a
\emph{simple component}~(\cite{P1}). This is, by definition, the
group $\Z$ with a positive cone $G^+$ such that $G=(\Z,G^+)$ is
partially ordered and simple. It was proved in~\cite[Proposition
2.4 (ii)]{RorVill} and~\cite[Proposition 2.5]{RorVill} that, if
$(\Z,G^+)$ is a simple component, then $G^+$ is the submonoid of
$\Z^+$ generated by a (unique and minimal) finite set of elements
$n_1,\ldots,n_k$ in $\Z^+$ (so that $G^+=\langle
n_1,\ldots,n_k\rangle$, and in fact $\gcd (n_1,\ldots,n_k)=1$). In
the particular case that $G^+=\langle k,l\rangle$ (and thus $k$
and $l$ are coprime integers), then one can determine the smallest
non-negative integer $N$ for which $N+p\in G^+$ for all $p\geq 0$,
but $N-1\notin G^+$ (see~\cite[Lemma 2.3]{RorVill}); namely,
$N=kl-k-l$.

\begin{prop}
\label{Bergman 1} Every simple ordered group of rank one $(G,G^+)$
is isomorphic (as an ordered group) to a direct limit of a
directed system $(G_n,f_{n,n+1})$, such that $G_n=(\Z,
G_n^+)$ is a simple component for every $n$ in $\N$ and
$f_{n,n+1}\colon G_n\to G_{n+1}$ is an order-embedding.
\end{prop}
\begin{proof}
Since $G$ is an abelian group of rank one, we can assume without
loss of generality that $1\in G$. Thus, by~\cite[Lemma 2.4]{P1},
there exists a (unique) generalized integer $\gn =\prod_{k\geq
1}{a_k}$ such that $G\cong\Z_{\gn}$. For each $n\geq 1$,
let $b_n=\prod_{k=1}^{n}{a_k}$, and define
$H_n=(1/b_n)\Z$. Notice that $H_n\subset H_{n+1}$ for each
$n\geq 1$, and also that $G=\bigcup_{n\geq 1}H_n$. Now, for every
$n\geq 1$, let $g_{n,n+1}\colon H_n\to H_{n+1}$ denote the
natural inclusion map, and define $H_n^+=G^+\cap H_n$.

We claim that $(H_n, H_n^+)$ is a simple group for each $n\geq 1$.
To check this, pick a non-zero element $x$ in $H_n^+$, and let
$y\in H_n$ be any element. Since $(G,G^+)$ is a simple group,
there exists $m$ in $\N$ such that $-mx\leq_G y\leq_G mx$. Thus,
$mx-y, y+mx\in G^+\cap H_n=H_n^+$, so that the previous inequality
also holds in $H_n$, as desired.

We claim now that the map $g_{n,n+1}\colon H_n\to H_{n+1}$ is
an order-embedding for every $n\geq 1$. By definition, it is a
positive one-to-one map. Now, let $x\in H_n$ be an element such
that $g_{n,n+1}(x)\in H_{n+1}^+=G^+\cap H_{n+1}$. Since
$x=g_{n,n+1}(x)\in G^+$ and $x\in H_n$, we conclude that $x\in
H_n^+$.

Finally, for each $n\geq 1$, let $f_n\colon\Z\to H_n$,
given by multiplication by $(1/b_n)$. Define
$G_n^+=b_nH_n^+\subseteq \Z$, and $G_n=(\Z,
G_n^+)$. Then, $f_n\colon G_n\to H_n$ is an order-isomorphism,
and hence the group $G_n$ is a simple component. Moreover, for
each $n\geq 1$, the map \[f_{n,n+1}=f_{n+1}^{-1}\circ
g_{n,n+1}\circ f_n\colon G_n\to G_{n+1}\] is an order-embedding.
Hence, for each $n\geq 1$ we get a commutative diagram
\[
\begin{CD}
G_n @>{f_{n,n+1}}>> G_{n+1}\\
@V{f_n}VV @VV{f_{n+1}}V\\
H_n @>{g_{n,n+1}}>> H_{n+1}\\
\end{CD}
\]
whence the maps $f_n$ induce an order-isomorphism $f$ from
$\varinjlim ((\Z, G_n^+),f_{n,n+1})$ onto $(G,G^+)$, as
wanted.
\end{proof}

As mentioned in the Introduction, one of the main objectives
in~\cite{P2} was to study the embedding of a certain class of
simple partially ordered groups of rank one into simple Riesz
groups of rank one. Such groups are parametrized by a triple
$(A,B,\mathcal{H})$, where $\mathcal{H}$ is a sequence of simple
groups (basically $\Z$ with different positive cones) and $A$ and
$B$ sequences of positive integers, all subject to certain axioms.
The proof of the key embedding result, established
in~\cite[Theorem 2.11]{P2}, is based on the fact that these groups
are isomorphic to a direct limit of an inductive system $((\Z,
G_n^+),f_{n,n+1})$ such that, for every $n$ in $\N$, the map
$f_{n,n+1}\colon (\Z, G_n^+)\to (\Z, G_{n+1}^+)$ is an
order-embedding given by multiplication by a non-negative integer
$a_n$ (where $A=(a_n)_{n\geq 1}$). Thus, in view of
Proposition~\ref{Bergman 1}, we can strengthen~\cite[Theorem
2.11]{P2} as follows:

\begin{theor}
\label{2.11 nou} Let $(G,G^+)$ be a simple ordered group of rank
one, and let $\gn$ be the generalized integer associated to $G$.
Given any infinite generalized integer $\gm$ coprime with $\gn$,
there exist a simple Riesz group of rank one $({\wt G}(\gm),{\wt
G}^+(\gm))$ and a positive morphism $$\tau\colon G\to {\wt
G}(\gm)$$ such that:
\begin{enumerate}
\item The group ${\wt G}(\gm)$ is isomorphic to $\Z_{\gn\cdot
\gm}$ (as abelian groups). \item The map $\tau $ is an
order-embedding.\hspace{\fill} $\Box$
\end{enumerate}
\end{theor}

The next result was also pointed out by G. Bergman.

\begin{lem}
\label{Bergman 2} Let $G_1=(\Q, G_1^+)$ and $G_2=(\Q, G_2^+)$ be
partially ordered abelian groups, and let $f\colon G_1\to G_2$
be a positive map. Then $f$ is an order-embedding if and only if
it is an isomorphism of ordered groups.
\end{lem}
\begin{proof}
Clearly, since $f$ is a group morphism from $\Q$ to $\Q$, it is
identically zero or an isomorphism of abelian groups.

Suppose that $f$ is an order-embedding, so that in particular it
is one-to-one. Hence the previous observation implies that $f$ is
an isomorphism. But then we also have
\[f(G_1^+)=G_2^+\cap f(G_1)=G_2^+\cap \Q=G_2^+\,,\]
so that it is an order-isomorphism. The converse is obvious.
\end{proof}

A first consequence of Theorem~\ref{2.11 nou} and
Lemma~\ref{Bergman 2} is the following characterization of
embeddability of simple ordered groups into simple Riesz groups.
This will be an important result in the sequel.
\begin{theor}
\label{Bergman def} An ordered simple group of rank one $(G,G^+)$
which is not a Riesz group can be embedded into a simple Riesz
group of rank one if and only if $G\not\cong\Q$.
\end{theor}
\begin{proof}
First, assume that $(G,G^+)$ is a simple ordered group of rank one
which is not a Riesz group, and suppose that $G\cong\Q$. Assume
that $(H,H^+)$ is a simple Riesz group of rank one and that
$f\colon G\to H$ is an order-embedding. Then $(H,H^+)$ is
order-isomorphic to a subgroup of $(\Q,\Q^+)$ and
the composition of the isomorphism $\Q\cong G$ with $f$ and
the embedding of $H$ into $\Q$ provides a non-zero morphism
from $\Q$ to $\Q$. Evidently this must be an
isomorphism, which implies that $f$ is surjective. But then
\[
f(G^+)=f(G)\cap H^+=H\cap H^+=H^+\,,
\]
that is, $f$ is an order-isomorphism. This implies that $G$ is a
Riesz group, a contradiction.

Conversely, suppose that $G\ncong \Q$. Since the generalized
integer associated to $\Q$ is $\gn (\Q)=\prod_{p_i\in
\mathbb{P}}p_i^{\infty }$, where $\mathbb{P}$ is the set of all
non-negative prime numbers, we conclude that for the generalized
integer associated to $G$, say $\gn (G)=\prod_{p_i\in
\mathbb{P}}p_i^{n_i}$, there exists at least one prime number
$p_k$ so that $n_k<\infty$. Now, multiplication by $p_k^{n_k}$
defines an order-isomorphism from $(G,G^+)$ onto $(p_k^{n_k}G,
p_k^{n_k}G^+)$. Notice that $\gn (p_k^{n_k}G)=\gn (G)/p_k^{n_k}$,
so that $\gn (p_k^{n_k}G)$ is coprime with $p_k$. Hence,
applying Theorem~\ref{2.11 nou} we get an order-embedding from
$(p_k^{n_k}G, p_k^{n_k}G^+)$ into a simple Riesz group of rank one
$(H,H^+)$. Thus, the composition of both maps gives us an
order-embedding from $(G,G^+)$ into $(H, H^+)$, as desired.
\end{proof}

\section{Intervals in basic building blocks}

This section, technical in nature, aims at the study of certain
simple groups of rank one. These will be used as our basic
building blocks in the subsequent sections, by connecting them
through order-embeddings and forming various inductive limits. We
shall focus on the construction of proper intervals in these
groups such that a certain multiple (that can be controlled)
equals the positive cone.

Let $G$ be a partially ordered abelian group with positive cone
$G^+$. A non-empty subset $X$ of $G^+$ is called an
\emph{interval} in $G^+$ if $X$ is upward directed and
order-hereditary.  We denote by $\Lambda(G^+)$ the set of
intervals in $G^+$. Note that $\Lambda (G^+)$ becomes an abelian
monoid with operation defined by $X+Y=\{ z\in G^+\mid z\leq x+y
\mbox{ for some } x\mbox{ in } X, y\mbox{ in }Y \}$. An interval
$X$ in $G^+$ is said to be \emph{generating} if every element of
$G^+$ is a sum of elements from $X$. We say that $X$ in $\Lambda
(G^+)$ is \emph{countably generated} provided that $X$ has a
countable cofinal subset (i.e. there is a sequence $(x_n)$ of
elements in $X$ such that for any $x$ in $X$ there exists $n$ in
$\N$ with $x\leq_G x_n$). Notice that, since any interval is
upward directed, if $(x_n)$ is a countable cofinal subset
generating an interval $X$, then we can choose a countable cofinal
subset $(y_n)$ generating $X$ with the property that $y_n\leq_G
y_{n+1}$ for all $n\geq 1$. We shall in this case use the notation
$X=\langle y_n\rangle$. We denote by $\Lambda _{\sigma }(G^+)$ the
set of all countably generated intervals in $G^+$.

\begin{defi}
\label{grup basic} {\rm Let $p$ and $q$ be positive integers such
that $1<q<p-q$ and that $\gcd(q,p)=1$. Denote by
$A=\langle q, p-q\rangle$ the submonoid of $\Z^+$ generated by $q$
and $p-q$. Let $s\in A\setminus\{0\}$ and take $r$ in $\Z^+$ be such that
$1<r<s-r$ and $\gcd(r, s)=1$. Denote by $B=\langle r,
s-r\rangle$.

Next, define $M$ to be the submonoid of $\Q^+$ whose
generators are fractions of the form
\[
\frac{k}{r}\,, \text{ where } k\in A\,, \text{ and }
\frac{k'}{r}\left(\frac{s}{r}\right)^n\,, \text{ where } k'\in B
\text{ and } n\geq 1\,.
\]
Define $(G,G^+)$ to be the group $G=M+(-M)$ with positive cone
$G^+=M$. Since $G$ is also directed, it follows that $G$ is a
simple partially ordered abelian group. Notice that $(G,G^+)$ is
not a Riesz group.}
\end{defi}

For all $n$ in $\N$, denote $e_n=(\frac{s}{r})^n$.

\begin{lem}
\label{interval} The set $D=\{x\in G^+\mid x\leq_G e_n \mbox{ for
some }n\}$ is a proper interval in $G^+$ such that $rD=G^+$.
\end{lem}
\begin{proof}
We first show that the sequence $(e_n)$ is increasing. By
construction, $e_n\in M$ for all $n$. Also, if $n\geq 1$, we have
\[
e_{n+1}-e_n=\left(\frac{s}{r}\right)^n\left(\frac{s}{r}-1\right)=
\left(\frac{s}{r}\right)^n\left(\frac{s-r}{r}\right)\,,
\]
which is an element of $M$ since $s-r\in B$. This proves that $D$
is an interval in $G^+$.

We now prove that $s\notin D$, while it is clear that
$s=\frac{sr}{r}\in M$. This will entail that $D$ is proper. In
order to achieve this, we proceed by induction. We evidently have
that $e_1-s=\frac{s(1-r)}{r}\notin M$ (because $1-r<0$). Assume, by way of
contradiction, that $s\not\leq_G e_{m-1}$ for some $m\geq 2$, and
that $s\leq_G e_m$. This means that we can find a natural number
$n$, elements $k_l$ in $B$ for $l=1,\ldots,n$, and an element $k$
in $A$ such that
\begin{equation} \label{laprimera}
\left(\frac{s}{r}\right)^m-s=e_m-s=\sum_{l=1}^n
\frac{k_l}{r}\left(\frac{s}{r}\right)^l+\frac{k}{r}\,.
\end{equation}
We can obviously choose $n$ above so that $k_n\neq 0$. Since
$k_n\in B$, we obtain that $k_n\geq r$. Therefore, substituting
$k_n$ by $r$ in~(\ref{laprimera}) we get the following bound:
\[
\left(\frac{s}{r}\right)^m>\left(\frac{s}{r}\right)^m-s\geq
\left(\frac{s}{r}\right)^n\,.
\]
This implies that $n<m$.

Now, the right-hand side of~(\ref{laprimera}) belongs to
$r^{-(n+1)}\Z^+$. Hence, after multiplying by $r^{n+1}$ we
get that
\[
s^mr^{n+1-m}-sr^{n+1}=r^{n+1}\left(\left(\frac{s}{r}\right)^m-s\right)\in\Z^+\,.
\]
Since $\gcd(r,s)=1$, the above implies that $m\leq n+1$.
Thus $m=n+1$.

We now claim that $r\nmid k_n$ and that $k_n<s$. Assume first that
$r\mid k_n$. Then the right hand side of~(\ref{laprimera}) would
belong to $r^{-n}\Z^+$. Hence
\[
\frac{s^{n+1}}{r}-r^ns=r^n\left(\left(\frac{s}{r}\right)^{n+1}-s\right)\in
\Z^+\,,
\]
contradicting the fact that $r$ and $s$ have no common factors.

Also, if $k_n\geq s$, then substituting $k_n$ by $s$ in~(\ref{laprimera}) we get
 the following bound
\[
\left(\frac{s}{r}\right)^{n+1}-s\geq
\left(\frac{s}{r}\right)^{n+1}\,,
\]
which is impossible. The claim is therefore established.

From our claim and the fact that $B=\langle r, s-r\rangle$, it
follows that $k_n=s-r$. Indeed, if we write $k_n=ar+b(s-r)$ for
some positive integers $a$ and $b$, we know that $b\neq0$ since
$r\nmid k_n$. If $a\neq 0$, then $s>k_n\geq r+s-r=s$, which is
impossible. Hence $a=0$. If now $b\geq 2$, then $s>k_n\geq
2(s-r)$, so $s-2r<0$, in contradiction of our election of $r$ and
$s$.

Finally
\begin{align*}
\left(\frac{s}{r}\right)^{n+1}-s & =\sum_{l=1}^n
\frac{k_l}{r}\left(\frac{s}{r}\right)^l+\frac{k}{r}\\
&=\sum_{l=1}^{n-1}
\frac{k_l}{r}\left(\frac{s}{r}\right)^l+\frac{s-r}{r}\left(\frac{s}{r}\right)^n+
\frac{k}{r}\\
&= \sum_{l=1}^{n-1} \frac{k_l}{r}\left(\frac{s}{r}\right)^l
\left(\frac{s}{r}\right)^{n+1}-\left(\frac{s}{r}\right)^n+\frac{k}{r}\,.
\end{align*}
This implies that
\[
e_n-s=\sum_{l=1}^{n-1}
\frac{k_l}{r}\left(\frac{s}{r}\right)^l+\frac{k}{r}\in M\,,
\]
which contradicts our inductive hypothesis since $n=m-1$.
Therefore, by induction, $s\not\leq_G e_m$ for all $m$ and so $D\neq
G^+$.

Next, we prove that $rD=M$. First, we claim that $2^ke_n\leq_G
re_{n+k-1}$ for all $n$ and all $k$. Indeed, if $k=1$, then
$2e_n\leq_G re_n$ for all $n$ (since $r\geq 2$). Now assume that,
for some $k\geq 2$, we have $2^ke_n\leq_G re_{n+k-1}$ for all $n$.
Then
\[
e_{n+k-1}\frac{r(s-2r)}{r}=r\left(\frac{s}{r}\right)^{n+k-1}
\left(\frac{s}{r}-2\right)=re_{n+k}-2re_{n+k-1}\leq_G
re_{n+k}-2^{k+1}e_n\,.
\]
Notice that by our choosing of $r$ and $s$, we have
$s-2r=s-r-r>0$. Therefore the element $e_{n+k-1}\frac{r(s-2r)}{r}$
belongs to $M$, and hence $re_{n+k}-2^{k+1}e_n\in G^+$. By
induction, the claim is proved.

Now take $e_1$, which belongs to $D$ and is non-zero. Since $G$ is
simple, $e_1$ is an order-unit. Given $x$ in $G^+$, there is then
a natural number $n$ such that $x\leq_G ne_1$. Choose $k$ such that
$n<2^k$. Hence, using the previous claim we conclude that $x\leq_G
ne_1\leq_G 2^ke_1\leq_G re_k$. This shows that $G^+\subseteq rD$.
Since the inclusion $rD\subseteq G^+$ is obvious, we get equality.
\end{proof}

\begin{prop}\label{jodiointerval}
Let $D$ be the interval defined in Lemma~\ref{interval}. Then, for any
$t\leq r-1$, we have $tD\neq G^+$, and $rD=G^+$.
\end{prop}
\begin{proof}
We have already checked in Lemma~\ref{interval} that $rD=G^+$. We
proceed by induction on $t$ to prove that $s\notin tD$ for any
$1\leq t\leq r-1$. The case $t=1$ is taken care of by the proof of
Lemma~\ref{interval}. Next, assume that, if $i<t$, we have
$s\not\leq_G ie_m$ for all $m$. We will prove that $s\not\leq_G
te_m$ for all $m$, using induction on $m$.

Since $te_1-s=t\frac{s}{r}-s=\frac{s(t-r)}{r}\notin M$, we see
that $s\not\leq_G te_1$.

Now, assume that $m\geq 2$ and that $s\not\leq_G te_j$ for all
$j<m$. By way of contradiction, assume that $te_m-s\in M$. This
means that we can find a natural number $n$, and elements $k_l$ in
$B$ for $l=1,\ldots,n$, $k$ in $A$, such that
\begin{equation}
\label{lasegunda} t\left(\frac{s}{r}\right)^m-s=\sum_{l=1}^n
\frac{k_l}{r}\left(\frac{s}{r}\right)^l+\frac{k}{r}\,.
\end{equation}
Since $k_n\in B$, we have that $k_n\geq r$. The right-hand side
of~(\ref{lasegunda}) belongs to $r^{-(n+1)}\Z^+$. Hence
\[
tr^{n+1-m}s^m-r^{n+1}s=r^{n+1}\left(t\left(\frac{s}{r}\right)^m-s\right)\in\Z^+\
.
\]
This, coupled with the assumptions that $t<r$ and
$\gcd(r,s)=1$, implies that $m\leq n+1$. We first deal
with the case $m=n+1$. From~(\ref{lasegunda}), we get
\begin{equation}
\label{latercera}
\left(\frac{s}{r}\right)^n\left(\frac{ts-k_n}{r}\right)=
\left(\frac{s}{r}\right)
^n\left(t\,\,\frac{s}{r}-\frac{k_n}{r}\right)=\sum_{l=1}^{n-1}
\frac{k_l}{r}\left(\frac{s}{r}\right)^l+\frac{k}{r}+s\,.
\end{equation}
Since the right-hand side of the above belongs to $r^{-n}\Z^+$ we
have that $s^n\left(\frac{ts-k_n}{r}\right)\in\Z^+$. As $r\nmid
s$, we conclude that $r\mid ts-k_n$. Write $ts-k_n=t'r$ for some
$t'$ in $\Z^+$. Now we have $t'r+k_n-tr-t(s-r)=0$. Adding
$(s-r)r-(s-r)-r$ to this equality, we get
\begin{equation}
\label{inter}
\begin{split}
(s-r)r-(s-r)-r & =t'r+k_n-tr-t(s-r)+(s-r)r-(s-r)-r\\
&=r(t'-t-1)+k_n+(s-r)(r-1-t)\,.
\end{split}
\end{equation}
By~\cite[Lemma 2.3]{RorVill} applied to the submonoid $B$,
$(s-r)r-(s-r)-r\notin B$. On the other hand, $r-1-t\geq 0$ and so
$(s-r)(r-1-t)\in B$. Since $r$, $k_n\in B$, it follows
from~(\ref{inter}) that $t'-t-1<0$, that is, $t'<t+1$.

We now substitute $ts-k_n=t'r$ in~(\ref{latercera}). We obtain
\[
\left(\frac{s}{r}\right)^nt'
=\left(\frac{s}{r}\right)^n\left(\frac{t'r}{r}\right
)=\sum_{l=1}^{n-1}
\frac{k_l}{r}\left(\frac{s}{r}\right)^l+\frac{k}{r}+s\,,
\]
whence $t'e_n-s=t'\left(\frac{s}{r}\right)^n-s\in M$, an absurdity
since $n=m-1$ and $t'<t$.

Next we deal with the case $n=m+a$ where $a\geq 0$. We rewrite~(\ref{lasegunda})
 as
\[
t\left(\frac{s}{r}\right)^m-s=\sum_{l=1}^{m-1}
\frac{k_l}{r}\left(\frac{s}{r}\right)^l+\sum_{l=m}^{m+a}
\frac{k_l}{r}\left(\frac{s}{r}\right)^l+\frac{k}{r}\,,
\]
that is,
\begin{equation}
\label{lacuarta}
\begin{split}
& \sum_{l=1}^{m-1}
\frac{k_l}{r}\left(\frac{s}{r}\right)^l+\frac{k}{r}+s\\
&=\left(\frac{s}{r}\right)^m\left(t-\sum_{l=m}^{m+a}
\frac{k_l}{r}\left(\frac{s}{r}\right)^{l-m}\right)\\
&=\left(\frac{s}{r}\right)^m\left(\frac{tr^{a+1}-\sum_{l=m}^{m+a}
k_ls^{l-m}r^{a+1-l+m}}{r^{a+1}}\right)\,.
\end{split}
\end{equation}
Let $u=\sum_{l=m}^{m+a} k_ls^{l-m}r^{a+1-l+m}$. Since the
left-hand side in~(\ref{lacuarta}) belongs to
$r^{-m}\Z^+$, we obtain (after multiplying the right-hand
side of the equality by $r^m$) that $r^{a+1}\mid tr^{a+1}-u$.
Write $tr^{a+1}-u=t''r^{a+1}$, for $t''$ in $\Z^+$, and rearrange as
$r^{a+1}(t''-t)+u=0$. Since $k_n=k_{m+a}\neq 0$, we have that
$u>0$. Therefore $t''<t$. Finally, we substitute
$tr^{a+1}-u=t''r^{a+1}$ in~(\ref{lacuarta}) and obtain
\[
\sum_{l=1}^{m-1}\frac{k_l}{r}\left(\frac{s}{r}\right)^l+\frac{k}{r}=
\left(\frac{s}{r}\right)^m\frac{t''r^{a+1}}{r^{a+1}}-s=t''e_m-s\,,
\]
so that $s\leq_G t''e_m$ and $t''<t$, a contradiction.
\end{proof}

Proposition~\ref{Bergman 1} allows us to write the group $(G,G^+)$
as an inductive limit of simple components and order-embeddings.
Below we present this representation in a way more related to the
construction and that will be used in the next section.

\begin{prop}
\label{erlimite} The group $(G,G^+)$ can be realized as a direct
limit $\varinjlim ((\Z,G_i^+),f_{i,i+1})$, where $(\Z,G_i^+)$ are
simple components and the maps $f_{i,i+1}\colon G_i\to G_{i+1}$ are
order-embeddings given by multiplication by $r$.
\end{prop}
\begin{proof}
Let $G_0^+=A$ and set $G_i^+=rG_{i-1}^++s^iB$ if $i\geq 1$. Since
$\gcd(r,s)=1$, the groups $(\Z,G_i^+)$ are simple components for all
$i$, and the maps $f_{i,i+1}$ given by multiplication by $r$ are
order-embeddings~(\cite[Lemma 2.3]{P2}).

Next, define $H_0=\frac{1}{r}\Z$, $H_0^+=\frac{1}{r}A$, and
$H_i=(\frac{1}{r})^i\Z$ and
$H_i^+=H_{i-1}^++\frac{1}{r}(\frac{s}{r})^iB$ if $i\geq 1$.

Clearly, we have the following commutative diagram:
\[
\begin{CD}
(H_0,H_0^+) @>>> (H_1,H_1^+)@>>> (H_2,H_2^+)@>>>\cdots\\
@V{r\cdot}VV  @VV{r^2\cdot}V @VV{r^3\cdot}V \\
(\Z,G_0^+) @>{f_{0,1}}>> (\Z,G_1^+)@>{f_{1,2}}>> (\Z,G_2^+)@>{f_{2,3}}>>\cdots\
\
\end{CD}
\]
where the maps in the top row are given by inclusions and all
columns are order-isomorphisms. The limit of the top row is
$(G,G^+)$ and it follows that the natural induced map to the limit
of the bottom row is an order-isomorphism, as desired.
\end{proof}

\begin{rema}
\label{clase} {\rm Following~\cite[Section~2]{P2}, it is easy to
see that $(G,G^+)$ is order-isomorphic to the group
$(G(A',B',\mathcal{H}'),G^+(A',B',\mathcal{H}'))$ associated to
the data triple $$(A',B',\mathcal{H}')=((r)_{i\geq 1},
(s^i)_{i\geq 1},\{A\}\cup \{B\}_{i\geq 2}).$$}
\end{rema}

\section{A first wild example}

In this section we construct our first example of a simple Riesz
group $(G,G^+)$ that contains an (even) decreasing sequence of
(unbounded) intervals $(D_n)$ such that, the larger $n$ is, the
more copies of $D_n$ we have to add in order to get $G^+$. The
main ingredient is the construction carried out in the previous
section, which is exploited with a certain recurrence using
(infinite) commutative diagrams.

\begin{lem}
\label{imageinterval} Let $f\colon (G,G^+)\to (H,H^+)$ be a
positive morphism. Let $D\subseteq G^+$ be an interval, and define
$D_f=\{x\in H^+\mid x\leq_H f(y)\mbox{ for some }y\in D\}$. Then:
\begin{enumerate}
\item $D_f$ is an interval. \item If $D$ is countably generated by
a sequence $(e_n)$, then $D_f$ is also countably generated, by
$(f(e_n))$. \item If $f$ is an order-embedding and $tD\neq G^+$
for some $t$ in $\N$, then $tD_f\neq H^+$. \item Let $r\in\N$.
Assume that $D$ is non-zero, $H$ is simple and $f$ is an
order-embedding. If $rD=G^+$ then $rD_f=H^+$.
\end{enumerate}
\end{lem}
\begin{proof}
(i) Evidently $D_f$ is non-empty. Let $x\in D_f$ and assume $0\leq
y\leq x$. By construction, there is an element $z$ in $D$ such
that $x\leq f(z)$, hence $y\leq f(z)$ and $y\in D_f$. This proves that
$D_f$ is order-hereditary. Next, take $x$, $y$ in $D_f$. There are
then $z_1$, $z_2$ in $D$ such that $x\leq f(z_1)$ and $y\leq
f(z_2)$. Since $D$ is an interval, there is a $z$ in $D$ such that
$z_i\leq z$. Hence $x,y\leq f(z)$, and $f(z)\in D_f$.

(ii) This is trivial.

(iii) Assume that $f$ is an order-embedding. Take $x$ in $G^+\setminus
tD$. Then $f(x)\in H^+\setminus tD_f$. For, if $f(x)\in tD_f$
there would be an element $y$ in $D_f$ such that $f(x)\leq ty$.
But then we could find an element $d$ in $D$ such that $y\leq
f(d)$, hence $f(x)\leq f(td)$. Since $f$ is an order-embedding,
this yields $x\leq td$, a contradiction.

(iv) Let $x\in D_f\setminus\{0\}$, and take $z_0$ in $D$ such that
$x\leq f(z_0)$. Since $H$ is simple, we know that $x$ is an
order-unit. If now $y\in H^+$, there is $n$ in $\N$ such that
$y\leq_H nx\leq_H nf(z_0)=f(nz_0)$. Now, $nz_0\in G^+=rD$, so that
we can find $z$ in $D$ for which $nz_0\leq rz$. Thus $y\leq_H
f(nz_0)\leq_H f(rz)=rf(z)$. Hence $y\in rD_f$.
\end{proof}
\begin{defidis}
\label{cons} {\rm Let $p$ and $q$ be positive integers such that
$1<q<p-q$ and that $\gcd(q,p)=1$. Set $A=\langle q, p-q\rangle$ as
in the Definition~\ref{grup basic}. Suppose that $(H,H^+)$ is a
simple ordered group of rank one such that there is an
order-embedding $(\Z, A)\hookrightarrow (H,H^+)$. By
Proposition~\ref{Bergman 1}, $(H,H^+)=\varinjlim ((\Z,G_{0,j}^+),
f_{0,j})$ with $(\Z, G_{0,0}^+)=(\Z, A)$ in such a way that
$(\Z,G_{0,j}^+)$ is a simple component and $f_{0,j}\colon
(\Z,G_{0,j}^+)\to (\Z,G_{0,j+1}^+)$ is an order-embedding, given
by multiplication by a non-negative integer $n_j$, for all $j\geq
0$. Let $\gn =\prod_{j\geq 0}n_j$ be the generalized integer
associated to the sequence $(n_j)$ and assume there exist a
positive integer $r>q$ such that $r$ is coprime with $\gn$, and a
positive integer $s$ in $A$ such that $\gcd(r,s)=1$ and $r<s-r$.
Put $B=\langle r, s-r\rangle$.

Next, define $G_{i,0}^+=rG_{i-1,0}^++s^iB$ for $i>0$, and
$G_{i,j}^+=rG_{i-1,j}^++n_{j-1}G_{i,j-1}^+$ for $i,j>0$. Let
$f_{i,j}\colon(\Z, G_{i,j}^+)\to(\Z,G_{i,j+1}^+)$ be the morphism given
by multiplication by $n_j$, and let $g_{i,j}\colon(\Z,
G_{i,j}^+)\to(\Z,G_{i+1,j}^+)$ be given by multiplication by $r$.
Denote by $G_i^+=G_{i,i}^+$, and
$f_i=g_{i,i+1}f_{i,i}=f_{i+1,i}g_{i,i}$\,.

Define $(K,K^+)=\varinjlim (\Z,G_{i,0}^+)$ and notice that, by
Remark~\ref{clase}, $(K,K^+)$ belongs to the class introduced
in~\cite[Definition 2.1]{P2}. It follows then by~\cite[Proposition
2.5]{P2} that $(K,K^+)$ is a simple group of rank one.
Observe that this construction yields the following commutative
diagram of groups and group morphisms:
\[
\begin{CD}
(\Z,G_{0,0}^+) @>{f_{0,0}}>> (\Z,G_{0,1}^+)@>{f_{0,1}}>> (\Z,G_{0,2}^+)@>{f_{0,2
}}>>\cdots\\
@V{g_{0,0}}VV  @VV{g_{0,1}}V @VV{g_{0,2}}V \\
(\Z,G_{1,0}^+) @>{f_{1,0}}>> (\Z,G_{1,1}^+)@>{f_{1,1}}>> (\Z,G_{1,2}^+)@>{f_{1,2
}}>>\cdots\\
@V{g_{1,0}}VV  @VV{g_{1,1}}V @VV{g_{1,2}}V \\
\vdots & & \vdots & & \vdots\\
\end{CD}
\]
}
\end{defidis}

\begin{prop}
\label{ladiagonal} For the construction in~\ref{cons}, the following
conditions hold:
\begin{enumerate}
\item $(\Z,G_{i,j}^+)$ is a simple component for all choices of
$i$ and $j$. \item The morphisms $f_{i,j}$ and $g_{i,j}$ are
order-embeddings for all choices of $i$ and $j$. \item Let
$(G,G^+)$ be the direct limit of the inductive system
$((\Z,G_i^+), f_i)$. Then $(G,G^+)$ is a simple group of rank one
and there are order-embeddings from $(H,H^+)$ into $(G,G^+)$ and
from $(K,K^+)$ into $(G,G^+)$. \item There exists an interval
$D_r\subseteq G^+$ such that $tD_r\neq G^+$ for $t\leq r-1$ and
$rD_r=G^+$.
\end{enumerate}
\end{prop}
\begin{proof}
(i) Since $\gcd(r,s)=1$,~\cite[Lemma 2.3 (1)]{P2} ensures that
$(\Z,G_{1,0}^+)$ with $G_{1,0}^+=rA+sB$ is a simple component.
Assume now that $(\Z,G_{i,0}^+)$ is a simple component. Since
$G_{i+1,0}^+=rG_{i,0}^++s^iB$, we can use~\cite[Lemma 2.3(1)]{P2}
again to conclude that $(\Z,G_{i+1,0}^+)$ is also a simple
component. Hence, by induction $(\Z,G_{i,0}^+)$ is a simple
component for all choices of $i$.

Next, assume that for all $k<j$, we have that $(\Z,G_{i,k}^+)$ is
a simple component for all $i$. We want to prove that $(\Z,
G_{i,j}^+)$ is a simple component for all $i$. By the discussion
in~\ref{cons}, we know that $(\Z, G_{0,j}^+)$ is a simple
component. Now, suppose that $(\Z,G_{i,j}^+)$ is a simple
component for some $i$. Then, since
$G_{i+1,j}^+=rG_{i,j}^++n_{j-1}G_{i+1,j-1}^+$ and
$\gcd(r,n_{j-1})=1$, another application of~\cite[Lemma
2.3(1)]{P2} allows us to conclude that $(\Z,G_{i+1,j}^+)$ is also
a simple component. The proof is then complete by induction.

(ii) By assumption, $f_{0,j}$ is an order embedding for all $j$.
Notice also that $g_{i,0}$ is also an order-embedding for all $i$,
by Proposition~\ref{erlimite}. Assume that $f_{i,j}$ is an order
embedding, and consider the following diagram:
\[\begin{CD} (\Z,G_{i,j}^+) @>{f_{i,j}=n_j\cdot}>> (\Z,G_{i,j+1}^+)\\
@V{g_{i,j}=r\cdot}VV  @VV{g_{i,j+1}=r\cdot}V\\
(\Z,G_{i+1,j}^+) @>{f_{i+1,j}=n_j\cdot}>> (\Z,G_{i+1,j+1}^+)
\end{CD}\]

Since
$G_{i+1,j+1}^+=f_{i+1,j}(G_{i+1,j}^+)+g_{i,j+1}(G_{i,j+1}^+)$, we
conclude from~\cite[Proposition 2.10]{P2} that $f_{i+1,j}$ and
$g_{i,j+1}$ are also an order-embeddings. Hence, it follows by
induction that $f_{i,j}$ and $g_{i,j}$ are order-embeddings for
all choices of $i$ and $j$.

(iii) That $(G,G^+)$ is a simple group follows from~\cite[Lemma
2.4]{P2}. By~\cite[Lemma 2.4]{P1}, $G$ is isomorphic to $\Z_{\gn
r^{\infty }}$, and so it is a group of rank one.

For every $j$, let $g_j\colon(\Z,G_{0,j}^+)\to (\Z,G_j^+)$ be defined
by $g_j=g_{j-1,j}\cdots g_{1,j}g_{0,j}$. Then $g_j$ is an
order-embedding and $g_{j+1}f_{0,j}=f_jg_j$ for all $j$. It
follows then from~\cite[Lemma 2.9]{P2} that the naturally induced
map $\varphi\colon(H,H^+)\to (G,G^+)$ is an order-embedding. In a
similar fashion, there is an order-embedding $\psi\colon(K,K^+)\to
(G,G^+)$.

(iv) By Proposition~\ref{erlimite}, together with
Lemma~\ref{interval} and Proposition~\ref{jodiointerval}, $K^+$
contains an interval $D$ such that $tD\neq K^+$ for $t\leq r-1$
and $rD=K^+$. Therefore, if we let $D_r=D_{\psi}$, then conditions
(iii) and (iv) in Lemma~\ref{imageinterval} ensure that $D_r$ will
do the job.
\end{proof}

\begin{theor}
\label{ertorema} Let $J=(q_i)_{i\geq 1}$ be a sequence of
non-negative, relatively prime integers. Let $I=(a_j)_{j\geq 1}$ be a
sequence such that every $a_j\in J$, while each $q_i$ in $J$ appears
infinitely many times in $I$. Let $\gn =\prod_{k\geq
1}q_k^{\infty }$ be the generalized integer associated to $J$.
Then, for any infinite generalized integer $\gm$ that is coprime
with $\gn$, there exists a simple Riesz group of rank one $G(\gm
)$ such that:
\begin{enumerate}
\item For every $q_{i}\in J$, there is a countably generated
interval $D_i$ satisfying $tD_i\neq G(\gm )^+$ for $t\leq q_i-1$
and $q_iD_i=G(\gm)^+$. \item The group $G(\gm )$ is isomorphic to
$\Z_{\gn\cdot\gm}$ (as abelian groups).
\end{enumerate}
\end{theor}
\begin{proof}
We first construct a simple ordered group of rank one $(G,G^+)$
satisfying condition (i) and such that $G\cong\Z_{\gn}$. To do so,
we proceed inductively.

Take $p_1$ such that $p_1>2q_1$ and $\gcd(p_1,q_1)=1$. Let
$A=\langle q_1, p_1-q_1\rangle$. Take $p_2 $ in $A$ such that
$p_2>2q_2$ and $\gcd(p_2,q_2)=1$. Let $B_1=\langle q_2,
p_2-q_2\rangle$. We construct the following commutative diagram of
groups and group morphisms:
\begin{equation}
\label{diagram1}
\begin{CD}
(\Z,G_{0,0}^{(1)+}) @>{q_1\cdot}>>
(\Z,G_{0,1}^{(1)+})@>{q_1\cdot}>>
(\Z,G_{0,2}^{(1)+})@>{q_1\cdot}>>\cdots\\ @V{q_2\cdot}VV
@VV{q_2\cdot}V @VV{q_2\cdot}V \\ (\Z,G_{1,0}^{(1)+})
@>{q_1\cdot}>> (\Z,G_{1,1}^{(1)+})@>{q_1\cdot}>>
(\Z,G_{1,2}^{(1)+})@>{q_1\cdot}>>\cdots\\ @V{q_2\cdot}VV
@VV{q_2\cdot}V @VV{q_2\cdot}V \\ \vdots  & & \vdots &  & \vdots \\
\end{CD}
\end{equation}
where $G_{0,0}^{(1)+}=A$, $G_{0,i}^{(1)+}=q_1
G_{0,i-1}^{(1)+}+p_1^iA$, $G_{i,0}^{(1)+}=q_2
G_{i-1,0}^{(1)+}+p_2^iB_1$. By~\cite[Lemma 2.3]{P2},
$(\Z,G_{0,i}^{(1)+})$ is a simple component for all $i$ and the
maps in the top row are order-embeddings. By
Proposition~\ref{ladiagonal}, all groups $(\Z,G_{i,j}^{(1)+})$ are
simple components and all maps in the diagram are
order-embeddings.

Let $(G_0,G_0^+)$ be the direct limit of the top row,
$(H_0,H_0^+)$ be the limit of the first column, and let
$(G_1,G_1^+)$ be the limit of the diagonal terms (under the
natural maps, obtained by composition). By~\cite[Proposition
2.5]{P2}, $(G_0,G_0^+)$ and $(H_0,H_0^+)$ are simple groups of
rank one. By condition (iii) in Proposition~\ref{ladiagonal}, $(G_1,
G_1^+)$ is also a simple group of rank one and there are
order-embeddings
\[
\tau_0\colon(G_0,G_0^+)\to (G_1,G_1^+)\,,\mbox{ and }
\psi_0\colon(H_0,H_0^+)\to (G_1,G_1^+)\,.
\]
By Lemma~\ref{interval}, Proposition~\ref{jodiointerval} and
Proposition~\ref{erlimite}, there are countably generated
intervals $D_1'$ in $G_0^+$ and $D^0$ in $H_0^+$ such that $tD_1'\neq
G_0^+$ if $t\leq q_1 -1$ and $q_1D_1'=G_0^+$; also, $tD^0\neq H_0^+$ if $t\leq q_2
 -1$, and $q_2D^0=H_0^+$.

Let $D_2'=(D^0)_{\psi_0}$. Then, Lemma~\ref{imageinterval}
ensures that $D_2'$ is a countably generated interval in $G_1^+$
such that $tD_2'\neq G_1^+$ for $t\leq
q_2-1$, and $q_2D_2'=G_1^+$.

Next, relabel the diagonal as the top row (i.e. let
$G_{0,i}^{(2)+}=G_{i,i}^{(1)+}$ for $i\geq 0$) and take $p_3$ in
$A$ such that $p_3>2q_3$ and $\gcd(p_3,q_3)=1$. Let
$B_2=\langle q_3, p_3-q_3\rangle$, and construct a commutative
diagram as before:
\begin{equation}
\label{diagram2}
\begin{CD}
(\Z,G_{0,0}^{(2)+}) @>{q_1q_2\cdot}>>
(\Z,G_{0,1}^{(2)+})@>{q_1q_2\cdot}>>
(\Z,G_{0,2}^{(2)+})@>{q_1q_2\cdot}>>\cdots\\ @V{q_3\cdot}VV
@VV{q_3\cdot}V @VV{q_3\cdot}V \\ (\Z,G_{1,0}^{(2)+})
@>{q_1q_2\cdot}>> (\Z,G_{1,1}^{(2)+})@>{q_1q_2\cdot}>>
(\Z,G_{1,2}^{(2)+})@>{q_1q_2\cdot}>>\cdots\\ @V{q_3\cdot}VV
@VV{q_3\cdot}V @VV{q_3\cdot}V \\ \vdots  & & \vdots &  & \vdots \\
\end{CD}
\end{equation}
Observe that, by construction, $(G_1,G_1^+)$ is the inductive
limit of the first row. Let $(H_1, H_1^+)$ be the inductive limit
of the first column and $(G_2,G_2^+)$ the inductive limit of the
diagonal terms. The same line of argument as before shows that
$(H_1,H_1^+)$ and $(G_2,G_2^+)$ are simple groups of rank one and
that there are order-embeddings
\[
\tau_1\colon(G_1,G_1^+)\to (G_2,G_2^+)\,,\mbox{ and }
\psi_1\colon(H_1,H_1^+)\to (G_2,G_2^+)\,.
\]
Another application of Lemma~\ref{interval}, Proposition~\ref{jodiointerval} and
 Proposition~\ref{erlimite} provides us
with a countably generated interval $D^1$ in $H_1^+$ such that
$tD^1\neq H_1^+$ if $t\leq q_3-1$ and $q_3D^1=H_1^+$. Let
$D_3'=(D^1)_{\psi_1}$. Then $D_3'$ is also a countably
generated interval (in $G_2^+$), by Lemma~\ref{imageinterval},
that satisfies $tD_3'\neq G_2^+$ for $t\leq q_3 -1$ and
$q_3D_3'=G_2^+$.

Continuing in this way, we get a sequence of simple groups of rank
one and order-embeddings
\[
\begin{CD}
(G_0,G_0^+) @>{\tau_0}>> (G_1,G_1^+)@>{\tau_1}>> (G_2,G_2^+)@>{\tau_2}>>\cdots\\
\end{CD}
\]
such that for each $i$, $G_i^+$ contains a countably generated
interval $D_{i+1}'$ such that $tD_{i+1}'\neq G_i^+$ for $t\leq
q_{i+1} -1$ and $ q_{i+1}D_{i+1}'=G_i^+$.

Let $(G,G^+)$ be the limit of this inductive system. Denote by
$\ol\tau_i\colon(G_i,G_i^+)\to (G,G^+)$ the natural maps. Now
define $D''_{i+1}=(D_{i+1}')_{\ol\tau_i}$. By
Lemma~\ref{imageinterval}, all the intervals $D''_j$ will satisfy
$tD''_j\neq G^+$ for all $t\leq q_j$, and $q_jD''_j=G^+$.
By~\cite[Lemma 2.4]{P2}, $(G,G^+)$ is a simple group, and since
$G\cong \Z_{\gn}$ by construction, it is a group of rank one.

Now, given any infinite generalized integer $\gm$ coprime with $\gn$,
Theorem~\ref{2.11 nou} ensures the existence of a simple Riesz
group of rank one $(G(\gm ), G(\gm )^+)$ and an order embedding
$\tau\colon (G,G^+)\to (G(\gm ), G(\gm )^+)$ such that $G(\gm
)$ is isomorphic to $\Z_{\gn\cdot\gm}$ (as abelian groups), thus
proving condition (ii). For each $i\geq 1$ define
$D_i=(D_i'')_{\tau }$. Then, by Lemma~\ref{imageinterval},
for every $q_i$ in $I$, $D_i$ satisfies that $tD_i\neq
G(\gm )^+$ for $t\leq q_i-1$ and $q_i D_i=G(\gm )^+$.
\end{proof}

Let $(G,u)$ be a partially ordered abelian group with order-unit.
We denote by $S(G,u)$ (or by $S_u$ if no confusion may arise) the
compact convex space of \emph{states} on $(G,u)$, that is, the set
of group morphisms $s\colon G\to \R$ such that $s(u)=1$. We use
$\mathrm{Aff}(S_u)^+$ to refer to the monoid of positive, affine
and continuous functions from $S_u$ to $\R^+$, and  $\phi_u\colon
G^+\to \mathrm{Aff}(S_u)^+$ stands for the natural evaluation map.
Let $\mathrm{LAff}_{\sigma}(S_u)^{++}$ be the monoid of strictly
positive, affine, lower semicontinuous functions from $G^+$ to
$\R^+$ that are pointwise suprema of increasing sequences of
functions from $\mathrm{Aff}(S_u)^+$.

If $D$ is a fixed interval in $\Lambda_{\sigma }(G^+)$, we denote
by $\Lambda _{\sigma ,D}(G^+)$ the submonoid of $\Lambda _{\sigma
}(G^+)$ whose elements are intervals $X$ in $\Lambda _{\sigma
}(G^+)$ such that $X\subseteq nD$ for some $n$ in $\N$, and we
denote by $W_{\sigma }^D(G^+)$ the submonoid of $\Lambda _{\sigma
,D}(G^+)$ whose elements are intervals $X$ in $\Lambda _{\sigma
,D}(G^+)$ such that there exists $Y$ in $\Lambda _{\sigma
,D}(G^+)$ with $X+Y=nD$ for some $n$ in $\N$. If now $D$ is a
countably generated interval in $G^+$ that is also generating, set
$d=\sup\phi_u(D)$, and define (see~\cite{Perera})
\[
W_{\sigma}^d(S_u)=\{f\in\mathrm{LAff}_{\sigma}(S_u)^{++}\mid f+g=nd
\mbox{ for some }g\mbox{ in }\mathrm{LAff}_{\sigma}(S_u)^{++}\mbox{
and }n\mbox{ in }\N\}\,.
\]
The disjoint union $G^+\sqcup W_{\sigma}^d(S_u)$ can be endowed with
a monoid structure by extending the natural operations and setting
$x+f=\phi_u(x)+f$ whenever $x\in G^+$ and $f\in W_{\sigma}(S_u)$.

Recall that an interval $X$ in $G^+$ is said to be \emph{soft}
(see, e.g.~\cite{G-H}) provided that for each $x$ in $X$, there is
an element $y$ in $X$ and a natural number $n$ such that
$(n+1)x\leq_G ny$. Observe that in case the interval $X$ satisfies
$rX=G^+$, then $X$ is soft. Indeed, if $x\in X\setminus\{0\}$,
then $(r+1)x\in G^+=rX$, hence there is an element $y$ in $X$ such
that $(r+1)x\leq_G ry$.

It was proved in~\cite[Theorem 3.8]{Perera} that, if $(G,u)$ is a
simple Riesz group with order-unit, and $D$ is a non-zero, soft,
countably generated interval in $G^+$, then the map
\begin{equation}
\label{lafi}
\varphi\colon W_{\sigma }^D(G^+)\to G^+\sqcup W_{\sigma
}^d(S_u)
\end{equation}
given by the rule $\varphi ([0,x])=x$ for any $x$ in $G^+$, and by
$\varphi (X)=\sup \phi _u(X)$ for any soft interval $X$ in
$W_{\sigma}^D(G^+)$, is a normalized monoid morphism. It becomes
an isomorphism if $G$ satisfies some additional assumptions,
namely if $G$ is non-atomic and strictly unperforated. Recently,
the first and second authors have shown that injectivity is
equivalent to strict unperforation~\cite[Theorem 3.2]{op}, and
surjectivity corresponds to a special property satisfied by the
generating interval $D$~\cite[Theorem 3.5]{op}. If $D$ is a soft
generating interval such that $\varphi(D)$ is identically
infinite, then we say that a soft interval $X$ in
$W_{\sigma}^D(G^+)$ is \emph{unbounded} provided that
$\varphi(X)=\sup\phi_u(X)=\infty$. Notice that this does not
depend on the choice of the order-unit. If $v$ is another
order-unit for $G$, then it follows from~\cite[Proposition
6.17]{POAG} that the state spaces $S_u$ and $S_v$ are
homeomorphic.

For the proof in the result below, we recall the following definition: An
abelian monoid $M$ is a \emph{refinement monoid} if, for all
$x_1,x_2,y_1,y_2$ in $M$ that satisfy $x_1+x_2=y_1+y_2$, there
exist elements $z_{ij}$ in $M$, for $i,j=1,2$, such that
$\sum\limits_{j=1}^2 z_{ij}=x_i$ and $\sum\limits_{i=1}^2
z_{ij}=y_j$ (see, e.g.~\cite{Whr}).
\begin{prop}
\label{intervalos} Let $(G,G^+)$ be a simple Riesz group, let
$I=(q_i)_{i\geq 1}$ be an increasing sequence of non-negative
integers such that $\gcd(q_{i},q_{j})=1$ for all different $i$ and
$j$. For every $q_{i}$ in $I$, assume that $D_i$ is a countably
generated interval in $W_{\sigma}^D(G^+)$ satisfying $tD_i\neq
G^+$ for $t\leq q_{i}-1$ and $q_i D_i=G^+$. Then there exists a
descending sequence of intervals $(X_i)$ such that $tX_i\neq G^+$
for $t\leq q_{i}-1$, and $(\prod _{j=1}^{i}q_j)X_i=G^+$.
\end{prop}
\begin{proof}
Let $M=\Lambda_{\sigma} (G^+)$ be the monoid of countably
generated intervals in $G^+$, with the algebraic ordering that we
shall denote by $\leq_M$\,. By~\cite[Proposition 2.5]{K0Good},
$\Lambda_{\sigma} (G^+)$ is a refinement monoid. Let $X_1=D_1$.
Since $q_2 D_2=G^+$, we have $X_1+G^+= q_2 D_2$. Hence,
by~\cite[Lemma 1.9]{Whr}, there exist intervals
$X_{11},X_{12},...,X_{1q_{2}}\leq_M D_2$ such that
$X_1=X_{11}+X_{12}+...+X_{1q_{2}}$ and $X_{11}\leq_M X_{12}\leq_M
\cdots \leq_M X_{1q_{2}}\leq_M X_1$. Let $X_2=X_{1q_2}$. Notice
that, if $tX_2=G^+$ for any $t\leq q_2-1$, then $tD_2=G^+$,
contradicting our assumption on $D_2$. Thus $tX_2\neq G^+$ for
$t\leq q_{2}-1$. Observe that
\[
X_1=X_{11}+\cdots+X_{1q_2}\leq_M X_{1q_2}+\cdots+X_{1q_2}=q_2X_2\,.
\]
Since $q_1X_1=q_1D_1=G^+$, we have that $(q_1q_2)X_2=G^+$. Now we
can apply the same argument to the equality $X_2+G^+ =q_3 D_3$, so
that we get an interval $X_3\leq_M X_2$ such that $tX_3\neq G^+$
for $t\leq q_{3}-1$ and $(q_1q_2q_3)X_3=G^+$. Continuing in this
way we get a descending sequence of intervals $(X_i)_{i\geq 1}$
such that $tX_i\neq G^+$ for $t\leq q_{i}-1$ and $(\prod
_{j=1}^{i}q_j)X_i=G^+$.
\end{proof}

\begin{rema}\label{kesinen} {\rm Notice that, if $(G,G^+,u)$ is a partially ordered abelian
group with order-unit, and $D\subseteq G^+$ is an interval such
that $nD=G^+$ for some natural number $n$, then $\varphi
(D)=\infty$, i.e. $D$ is an unbounded interval. To see this,
notice that, given any non-zero element $x$ in $G^+$, there exists
an element $y$ in $D$ such that $x\leq_G ny$. Hence, for any state $s$
on $G$ we have $0\leq s(x)\leq ns(y)$, i.e. $0\leq \frac{\phi
_u(x)}{n}\leq\phi _u(y)$. Thus, in order to see that $\varphi
(D)=\infty$, it is enough to show that $\varphi (G^+)=\infty$. But
now, for every $m$ in $\N$, we have that $mu\in G^+$, and then
$0<m=\phi _u(mu)$, whence $\varphi (G^+)=\infty$. In particular,
this fact applies to the intervals $D_i$, $X_j$ in Proposition~\ref{intervalos}.}
\end{rema}

The construction just carried out in Theorem~\ref{ertorema}
guarantees that we are in position to apply
Proposition~\ref{intervalos} and obtain a somewhat more refined
example, as follows.

\begin{theor}
\label{eldarrer} Let $J=(q_i)_{i\geq 1}$ be a sequence of
non-negative integers such that $\gcd(q_{i},q_{j})=1$ for all
$i,j\geq 1$ (such that $i\neq j$), and let $I=(a_j)_{j\geq 1}$ be
a sequence such that every $a_j\in J$, while each $q_i\in J$
appears an infinite number of times in $I$. Let $\gm$ be a
generalized integer, that is coprime with $\gn(I)$. Let
$(G,G^+)=(G(\gm),G^+(\gm))$ be the simple Riesz group of rank one
constructed in Theorem~\ref{ertorema}. Then $G$ contains a
descending sequence of unbounded intervals $(D_i)$ such that
$(\prod _{j=1}^{i}q_j)D_i=G^+$ for all $i$, while $tD_i\neq G^+$
whenever $t\leq q_i-1$.
\end{theor}
\begin{proof}
We only need to check that $(G,G^+)$ fulfills the hypotheses of
Proposition~\ref{intervalos}. The sequence of intervals obtained
in the conclusion of Theorem~\ref{ertorema}, say $(D_i')$,
satisfies that $tD_i'\neq G^+$ for all $t\leq q_i-1$, and
$q_iD_i'=G^+$ for all $i$. Thus, result holds by
Proposition~\ref{intervalos}.
\end{proof}

The examples we have just obtained could be considered as an
intermediate step towards constructing a simple Riesz group $(G,
G^+)$, together with an unbounded interval $D$ in $G^+$ such that
$nD\ne G^+$ for every $n$ in $\N$. In fact, a natural candidate
for such an interval could the intersection of the descending
chain of intervals appearing in Theorem~\ref{eldarrer}.
Unfortunately, even under the hypotheses of
Theorem~\ref{eldarrer}, we are not able to prove whether or not
the intersection $D=\bigcap _{i\geq 1}D_i$ is an interval or even
an unbounded subset of $G^+$, where $(D_i)$ is a descending
sequence of countably generated, unbounded soft intervals $(D_i)$
such that $tD_i\neq G^+$ for $t\leq q_i-1$ and $(\prod
_{j=1}^{i}q_j) D_i=G^+$ ( for every $i\geq 1$).

\section{Taylor-made gaps in simple components under order-embeddings}

In order to obtain our desired example (announced in the
Introduction) of a simple Riesz group of rank one $(G,G^+)$,
together with an unbounded (countably generated) interval $D$ in
$G^+$ such that $nD\neq G^+$ for all $n$, we adopt the basic
philosophy of~\cite[Section~3]{P1}. This consists of reducing the
essential properties that an interval should have to a finite set
of properties occurring in simple components. For this, we need to
have some control over those non-negative integers in a simple
component that its positive cone may contain, and we also need to
construct order-embeddings among simple components under which
this control is preserved.

In view of the considerations, made at the beginning of Section~1
and related to results on simple components (see~\cite{P1} and~\cite{RorVill}),
a possible way of getting the desired control is
to consider submonoids of the non-negative integers generated by
coprime positive integers, and order-embeddings among simple
components whose positive cones have this particular form, using~\cite[Lemma 2.3
(2)]{P2}. The basic idea consists of strengthening
some arithmetic properties in order to force the expression of
non-negative integers to become positive elements in a certain
simple component.

\begin{lem}
\label{lema1} Let $N\in \N$. Given $a$ in $\N$, there exist $p$,
$c$ and $d$ in $\N$ such that $\gcd(a,p)=\gcd(a,c)=\gcd(c,d)=1$,
$pc\equiv pd\equiv 1\, (\mathrm{mod}\,\, a)$, $p>N$, $pc>aN$ and
$d>\max \{(a-1)pc+a(N-1),ac\}$.
\end{lem}
\begin{proof}
Throughout the proof, denote by $\overline{x}$ the class of an
element in $\Z/n\Z$ for any $n$. For $p$ and $c$ in $\N$, it is clear that
$\gcd(p,a)=\gcd(c,a)=1 $ is equivalent to the fact that
$\overline{p}$ and $\overline{c}$ are invertible in $\Z/a\Z$.
Therefore, if we take $p$ and $c$ in such a way that
$\overline{p}=\overline{c}^{-1}\in\Z/a\Z$, we will have that
$pc\equiv 1\,(\mathrm{mod}\,\, a)$. It is clear that there exist
infinitely many numbers $p$ and $c$ satisfying the above. We can
then take $p$, $c>N$ and also $pc>aN$. By a similar line of
argument, once $p$ is fixed, there are infinitely many $d$ in $\N$
such that $pd\equiv 1\, (\mathrm{mod}\,\, a)$. For any one of this
choices we have that $\overline{d}=\overline{p}^
{-1}=\overline{c}$ in $\Z/a\Z$, whence $d-c$ is divisible by $a$,
that is $d=c+a k$ for some $k$ in $\N$. Now, in $\Z/c\Z$, this
says
$\overline{d}=\overline{c}+\overline{a}\overline{k}=\overline{a}\overline{k}$.
We can choose $k$ big enough such that $d>\max
\{(a-1)pc+a(N-1),ac\}$ and $\gcd(c,k)=1$. This will also guarantee
that $d$ is invertible in $\Z/c\Z)$, that is $\gcd(c,d)=1$.
\end{proof}

\begin{nota}
\rm{Let $(\Z,H^+)$ be a simple component. There is then a
(uniquely determined) element $N$ in $H^+$ such that $N-1\notin
H^+$, and $N+k\in H^+$ for all $k$ in $\Z^+$. We shall denote this
element by $N_H$.

For the rest of this section, let us fix a simple component
$(\Z,H^+)$. Given $a$ in $\N$ we can choose by Lemma~\ref{lema1}
elements $p$, $c$ and $d$ in $\N$ such that
$\gcd(a,p)=\gcd(a,c)=\gcd(c,d)=1$, $pc\equiv pd\equiv 1\,
(\mathrm{mod}\,\, a)$, $p\in H^+$, $pc>aN_H$ and
$d>\max\{(a-1)pc+a(N_H-1),ac\}$.

Let $G^+=aH^++p\langle c,d\rangle$. Since $\gcd(c,d)=1$, we have
that $( \Z,\langle c,d\rangle)$ is a simple component. Since
$\gcd(a,p)=1$, we can use~\cite[Lemma 2.3]{P2} and conclude that
$(\Z,G^+)$ is a simple component and the map $(\Z,H^+)\to
(\Z,G^+)$ defined by multiplication by $a$ is an order-embedding.}
\end{nota}

We shall use these notations in the Proposition below and in the
next section.

\begin{prop}
\label{boquetes} Let $i$ in $\Z$ be such that $0\leq i\leq a-1$
and let $x\notin H^+$. Then $ipc+ax\notin G^+$. In particular, if
we denote $L_H=\{l_0,l_1,\ldots l_{a-1}\}$ where
$l_i=ipc+a(N_H-1)$, it follows that $L_H \cap G^+=\emptyset$.
Moreover, all integers that are congruent to $i$
$(\mathrm{mod}\,\, a)$ and bigger than $l_i$ belong to $G^+$.
\end{prop}
\begin{proof}
Since multiplication by $a$ is an order-embedding, we see that
$at\notin G^+$ if and only if $t\notin H^+$. In particular,
$a(N_H-1)\notin G^+$. Moreover, any multiple of $a$ which is
bigger than $a(N_H-1)$ will belong to $G^+$, as it will have the
form $at$ with $t=N_H+k$ ($k$ in $\Z^+$), and so $t\in H^+$.

Assume now that we have $0\leq i\leq a-1$ and $x\notin H^+$. We
have to prove that $ipc+ax\notin G^+$. By way of contradiction, if
$ipc+ax\in G^+$, we would then have that $ipc+ax=ay+pz$ where
$y\in H^+$ and $z\in\langle c,d\rangle$, that is, $z=cz_1+dz_2$
for some positive integers $z_1$, $z_2$. We therefore have
\begin{equation}
\label{unaecuacion} ipc+ax=ay+pcz_1+pdz_2\,.
\end{equation}
We now claim that $y<x$. We already know that $y\neq x$ because
$y\in H^+$. Assume that $y\geq x+1$, so that $y=x+k$ with $k\geq
1$. We would then have that $ipc+ax=
ay+pcz_1+pdz_2=ax+ak+pcz_1+pdz_2$, whence $ipc=ak+pcz_1+pdz_2$.
Thus (since also $c<d$),
\[
0<a\leq ak=ipc-pcz_1-pdz_2\leq ipc-pcz_1-pcz_2=pc(i-(z_1+z_2))\,,
\]
and so $i-(z_1+z_2)\geq 0$. Notice that also
\[
0\equiv ak=ipc-(pcz_1+pdz_2)\equiv i-(z_1+z_2)\, (\mathrm{mod}\,\,
a)\,.
\]
This implies that $i-(z_1+z_2)=ar$ for some positive integer $r$.
But since $i\leq a-1$ by assumption we conclude that $r=0$, hence
$i=z_1+z_2\leq a-1$. Then
$ipc=ak+pcz_1+pdz_2>ak+pc(z_1+z_2)=ak+pci$, and so $0>ak>0$. This
contradiction establishes the claim.

Going back to equation~(\ref{unaecuacion}), we find that $0\leq
a(x-y)=pcz_1+pdz _2-ipc$. Since $x\notin H^+$, we have that
$y<x\leq N_H-1$. If $z_2\neq 0$, then
 $az_2-(a-1)\geq 1$ and (using that $d>ac$ and that $pc>aN_H$), we get
\begin{align*}
a((N_H-1)-y) & \geq a(x-y)=pcz_1+pdz_2-ipc\\
& \geq pdz_2-ipc>pacz_2-ipc\geq pacz_2-(a-1)pc\\
& =pc(az_2-(a-1))>aN_H(az_2-(a-1))>aN_H\,,
\end{align*}
which is clearly not possible.

It follows then that $y_2=0$. This means that
$a(x-y)=pcz_1-ipc=pc(z_1-i)$, whence $z_1>i$ (because $x>y$). But
then $a(x-y)=pc(z_1-i)\geq pc>aN_H$. This implies that $x-y\in
H^+$ and since $y\in H^+$, it follows that $x\in H^+$, contrary to
our assumption.

We have proved that $ipc+ax\notin G^+$, whenever $0\leq i\leq a-1$
and $x\notin H^+$. In particular, since $N_H-1\not\in H^+$, we
have that $l_i=ipc+a(N_H-1)\not\in G^+$, hence $L_H\cap G
=\emptyset$.

Let now $t$ in $\Z^+$ be such that $t\equiv i\,
(\mathrm{mod}\,\,a)$ and $t>ipc+a (N_H-1)$. Then, since $pc\equiv
1\,(\mathrm{mod}\,\,a)$, we have that $t-ipc\equiv 0\,
(\mathrm{mod}\,\, a)$ and so $a(N_H-1)<t-ipc=as$ for some $s$ in
$\Z^+$. Then $N_H-1<s$, and therefore $s\in H^+$ and $t=as+ipc\in
aH^++p\langle c,d\rangle=G^+$. It follows from this that any
integer congruent to $i$ $(\mathrm{mod}\,\, a)$ which is bigger
than $l_i$ can be written as $ipc+as$ where $s>N_H-1$ and so they
all belong to $G ^+$.
\end{proof}

\begin{corol}
\label{elultimo} Under the hypotheses and notation of
Proposition~\ref{boquetes}, we have $N_G=l_{a-1}+1$.
\end{corol}
\begin{proof}
It is clear that $l_{a-1}$ is the largest element of the set $L_H$
defined in Proposition~\ref{boquetes}. Let $x>l_{a-1}$. We
obviously have that $x\equiv k\, (\mathrm{mod}\,\, a)$, for some
$0\leq k\leq a-1$. Since $x>l_{a-1}>l_k$ we obtain (using
Proposition~\ref{boquetes}) that $x\in G^+$.
\end{proof}

\section{A new wild example}

The main objective of this section is to construct a simple Riesz
group $(G,G^+) $ of rank one such that its positive cone contains
an unbounded interval $D$ that satisfies $nD\neq G^+$ for all $n$
in $\N$. This will be done inductively by constructing a sequence
of simple components connected by order-embeddings. We first
establish a Lemma that will provide the inductive step in the
Theorem below. Given a simple component $(\Z,H^+)$, retain from
the previous section the notation $N_H$ for the (unique) element
in $H^+$ such that $N_H-1\notin H^+$ but $N_H +k\in H^+$ for all
positive integers $k$.

\begin{lem}
\label{unpaso} Let $(\Z,H^+)$ be a simple component, let $x_1$,
$y_1\in H^+$ be such that $y_1=x_1+1$, and let $a>N_H$. There
exists then a simple component $(\Z,G^+)$ satisfying:
\begin{enumerate}
\item $a\cdot\,\colon (\Z,H^+)\to (\Z,G^+)$ is an order-embedding
and $a^ 2N_H<N_G$. \item There is an element $y_2$ in $G^+$ such
that:
\begin{enumerate}
\item $y_2-1\in G^+$. \item $ay_1<_G y_2$ and $y_2-ay_1>aN_H$.
\item $(N_H-1)ax_1\nless_G (N_H-1)y_2$.
\end{enumerate}
\end{enumerate}
\end{lem}
\begin{proof}
Notice that $(N_H-1)y_1-(N_H-1)x_1=(N_H-1)\notin H^+$, whence
$(N_H-1)x_1\nleq_H (N_H-1)y_1$.

Choose $p$, $c$ and $d$ as in Lemma~\ref{lema1}. Letting
$G^+=aH^++p\langle c,d\rangle$, we have that $(\Z, G^+)$ is a
simple component and multiplication by $a$ is an order-embedding.

Write $L_H=\{l_0,l_1,\ldots,l_{a-1}\}$ as in
Proposition~\ref{boquetes}, so we have that any integer congruent
to $i$ $(\mathrm{mod}\,\, a)$ which is larger than $l_i$ belongs
to $G^+$. Note also that $N_G=l_{a-1}+1$, by
Corollary~\ref{elultimo}. This equals to $N_G=(a-1)pc+a(N_H-1)+1$,
and hence we have
\begin{align*}
N_G & =(a-1)pc+a(N_H-1)+1>(a-1)(pc+(N_H-1)) \\
& >(a-1)(aN_H+(N_H-1))>(a-1)(a+1)(N_H-1)\\ &=(a^2-1)(N_H-1)>a^2
N_H\,,
\end{align*}
proving condition (i).

Let $y_2=pc+ay_1$, and observe that $y_2\in G^+$, because $y_1\in
H^+$ by assumption. Notice also that $y_2-ay_1=pc>aN_H$, by the
election of $p$ and $c$. Since $pc\in G^+$, we see that
$ay_1<_Gy_2$, thus verifying condition (ii)(b).

Since $a>N_H$, it follows from Proposition~\ref{boquetes} that
$l_{N_H-1}\notin G^ +$. Therefore, the fact that
\[
(N_H-1)y_2-(N_H-1)ax_1=(N_H-1)pc+a(N_H-1)=l_{N_H-1}\notin G^+
\]
implies that $(N_H-1)ax_1\nleq_G (N_H-1)y_2$. Hence condition
(ii)(c) also holds.

It remains to verify condition (ii) (a). Since
$y_2-1=pc+ay_1-1\equiv 0\, (\mathrm{mod}\,\, a)$ and
$y_2-1=pc+ay_1-1\geq pc>aN_H>a(N_H-1)=l_0$,
Proposition~\ref{boquetes} ensures that $y_2-1\in G^+$, as
desired.
\end{proof}
\begin{theor}\label{ertorema2}
Let $A$ be a strictly ascending sequence of non-negative integers.
Then, for any generalized integer $\gm$ coprime with $\gn (A)$,
there exists a simple Riesz group of rank one $G(\gm)$ such that:
\begin{enumerate}
\item There is an unbounded countably generated interval $D$
satisfying $nD\neq G(\gm)^+$ for all $n$ in $\N$\,. \item For some
generalized integer $\gn$ dividing $\gn (A)$, the group $G(\gm)$
is isomorphic to $\Z_{\gn \gm}$ (as abelian groups).
\end{enumerate}
\end{theor}
\begin{proof}
First, we will inductively construct a sequence of simple
components and order-embeddings
\[
a_i\cdot\, \colon(\Z,H_i^+)\to (\Z,H_{i+1}^+)\,,
\]
together with a sequence $(y_i)$ in $\Z^+$ ($i\geq 1$) such that
\begin{itemize}
    \item[(a)] For all $i\geq 1$, $a_i\in A$\,.
    \item[(b)] For all $i\geq 1$, $a_i>N_{H_i
}>a_{i-1}^2N_{H_{i-1}}>(a_1^2)^{i-1}N_{H_1}$ for all $i\geq 1$.
Also $y_i\in H_i^+$, and the element $x_i=y_i-1\in H_i^+$ for all
$i$\,.
    \item[(c)] $(N_{H_i}-1)x_i\nless_{H_i} (N_{H_i}-1)y_i$\,.
    \item[(d)] $a_i y_i<_{H_{i+1}}y_{i+1}$\,.
    \item[(e)] $(N_{H_j}-1)a_{i-1}a_{i-2}\cdots a_jx_j\nleq_{H_i} (N_{H_j}-1)y_i
$ for all $j\leq i-1$\,.
\end{itemize}
Let $(\Z,H_1^+)$ be any simple component such that $1\notin
H_1^+$. Let $y_1$ in $H_1^+$ be such that $x_1=y_1-1\in H_1^+$
(for example, $y_1=N_{H_1}+1$). Choose $a_1$ in $A$ with
$a_1>\max\{N_{H_1},3\}$. Then Lemma~\ref{unpaso} provides us with
a simple component $(\Z,H_2^+)$ (where $N_{H_2}>a_1^2N_{H_1}$) and
an element $y_2\in H_2^+$ such that the element $x_2=y_2-1\in
H_2^+$, $a_1y_1<_{H_2}y_2$, $(N_{H_1}-1)a_1x_1\nless_{H_1}
(N_{H_1}-1)y_2$ and $a_1N_{H_1}<y_2-a_1y_1$. Moreover,
multiplication by $a_1$ is an order-embedding from $(\Z, H_1^+)$
into $(\Z, H_2^+)$.

Suppose that $a_1,\ldots,a_{n-1}$, $H_1^+,\ldots,H_n^+$ and
$y_1,\ldots,y_n$ have been constructed satisfying conditions
(a)-(e) above.

Choose $a_n$ in $A$ with $a_n>N_{H_n}$, and apply
Lemma~\ref{unpaso} to obtain an order-embedding
\[
a_n\cdot\colon(\Z,H_n^+)\to (\Z,H_{n+1}^+)\,,
\]
where $(\Z,H_{n+1}^+)$ is a simple component, such that
$N_{H_{n+1}}>a_n^2N_{H_n}$. Moreover, there is an element
$y_{n+1}$ in $H_{n+1}^+$ such that the element
$x_{n+1}=y_{n+1}-1\in H_{n+1}^+$, $a_ny_n<_{H_{n+1}}y_{n+1}$,
$y_{n+1}-a_ny_n>a_nN_{H_n}$ and
$(N_{H_n}-1)a_nx_n\nless_{H_{n+1}}(N_{H_n}-1)y_{n+1}$. Hence
conditions (a)-(d) are satisfied (as well as condition (e) with
$i=n+1$ and $j=n$).

Thus, in order to see that condition (e) also holds with $i=n+1$,
we only need to consider the cases where $j\leq n-1$.

Notice that $y_{n+1}=p_nc_n+a_ny_n$ (by the proof of
Lemma~\ref{unpaso}) where $p_n$, $c_n$ are chosen as in
Lemma~\ref{lema1}. By our induction hypothesis,
\[
(N_{H_j}-1)y_n-(N_{H_j}-1)\prod\limits_{k=j}^{n-1}a_kx_j\notin
H_n^+\,,
\]
whenever $j\leq n-1$.

Since $N_{H_j}-1<a_n-1$ if $j\leq n-1$, Proposition~\ref{boquetes}
applies and so
\[
(N_{H_j}-1)p_nc_n+a_n[(N_{H_j}-1)y_n-(N_{H_j}-1)\prod\limits_{k=j}^{n-1}a_kx_j]\notin
H_{n+1}^+\,,
\]
that is,
\[
(N_{H_j}-1)y_{n+1}-(N_{H_j}-1)\prod\limits_{k=j}^na_kx_j\not\in
H_{n+1}^+
\]
for every $j\leq n-1$, as desired.

Next, let $(G,G^+)=\varinjlim ((\Z,H_i^+),a_i\,\cdot)$, and denote
by $f_n\colon (\Z, H_n^+)\to (G,G^+)$ the natural
(order-embedding) maps. By condition (d), $y_{i+1}-a_iy_i\in
H_i^+\setminus\{0\}$, hence
$f_{i+1}(y_{i+1})-f_i(y_i)=f_{i+1}(y_{i+1}-a_iy_i)\in
G^+\setminus\{0\}$. This shows that the interval $E=\langle
f_i(y_i)\rangle$ is soft and countably generated (see,
e.g.~\cite[Lemma 3.4]{Perera}).

Let $u=f_1(y_1)$ in $G^+$, and take this as an order-unit. Denote
by $s$ the (unique) state on $(G,u)$; for $i$ in $\N$, let $s_i$
denote the unique state on the simple component $(\Z, H_i^+)$ with
respect to the order-unit $u_i=a_{i-1}a_{i-2}\cdots a_1y_1$. We
now check that $E$ is unbounded, that is, $\sup\phi_u(E)=\infty$.
By the first part of the proof, $y_{i+1}=p_ic_i+a_iy_i$ where
$p_i$ and $c_i$ are chosen in such a way that $p_ic_i>a_iN_{H_i}$.
Then, by using condition (b) recurringly, we get
\begin{align*}
s(f_{i+1}(y_{i+1}))= &\, s_{i+1}(y_{i+1})  =\,
\frac{p_ic_i}{a_ia_{i-1}\cdots
a_1y_1}+\frac{a_iy_i}{a_ia_{i-1}\cdots a_1y_1}  >
\frac{p_ic_i}{a_ia_{i-1}\cdots a_1y_1} \\ > &
\frac{a_iN_{H_i}}{a_ia_{i-1}\cdots
a_1y_1}=\frac{N_{H_i}}{a_{i-1}\cdots a_1y_1} >
\frac{a_{i-1}^2N_{H_{i-1}}}{a_{i-1}\cdots
a_1y_1}=\frac{a_{i-1}N_{H_{i-1}}}{a_{i-2}\cdots
a_1y_1} \\ > & \frac{a_{i-1}a_{i-2}^2N_{H_{i-2}}}{a_{i-2}\cdots a_1y_1}
>  \cdots>\frac{a_{i-1}\ldots a_2a_1N_{H_1}}{y_1}>a_1^{(i-2)(i-1)}\frac{N_{H_1}}
{y_1}\,,
\end{align*}
and so clearly $\sup\phi_u(E)=\infty$.

Now, suppose that $nE=G^+$ for some $n$. Choose $j$ in $\N$ such
that $(N_{H_j}-1)\geq n$. We have that
$f_j((N_{H_j}-1)x_j)<_G(N_{H_j}-1)f_i(y_i)$ for all (suitably)
large $i$. This will happen in particular for some $i>j$, which
translates into $f_i((N_{H_j}-1)a_{i-1}a_{i-2}\cdots a_jx_j)<_G
f_i((N_{H_j}-1)y_i)$. Since $f_i$ is an order-embedding, we get
$(N_{H_j}-1)a_{i-1}a_{i-2}\cdots a_jx_j<_{H_i} (N_{H_j}-1)y_i$ for
some $i>j$, in contradiction with condition (e).

Hence, we have constructed a simple group of rank one $(G,G^+)$,
containing an interval $E\subset G^+$ such that $\varphi
(E)=\infty$ and $nE\ne G^+$ for every $n$ in $\N$. Notice that
$A'=(a_i)_{i\geq 1}$ is an infinite subsequence of $A$, so that
$\gn =\gn(A')$ is a generalized integer dividing $\gn (A)$.
Moreover, by construction, $G\cong \Z_{\gn }$ (as abelian groups).
Thus, for any generalized integer $\gm$ coprime with $\gn $, there
exists by Theorem~\ref{2.11 nou} a simple Riesz group of rank one
$(G(\gm),G^+(\gm))$ such that $G(\gm )\cong \Z _{\gn \gm}$, and an
order-embedding $\tau\colon G\to G(\gm)$. Then, by condition (iii)
in Lemma~\ref{imageinterval}, the interval $D=E_{\tau }=\langle
(\tau f_i)(y_i)\rangle$ satisfies that $nD\ne G^+(\gm)$ for every
$n$ in $\N$. Let $u$ be an order-unit in $G$. Since both $S(G, u)$
and $S(G(\gm), \tau (u))$ are singletons with (unique) states $s$
and $s'$ respectively, the affine continuous map $S(\tau )\colon S(G,
u)\to S(G(\gm), \tau (u))$ is an homeomorphism with
$S(\tau )(s')=s'\tau =s$. Hence, $\sup s'((\tau
f_i)(y_i))=\sup(s'\tau) (f_i(y_i))=\sup s(f_i(y_i))=\infty$,
whence $D$ is also unbounded. This completes the proof.
\end{proof}

\section{The monster example}

In this section, we will use the constructions carried out in
Theorems~\ref{ertorema} and~\ref{ertorema2} in order to construct
an example of a simple Riesz group of rank one containing
unbounded intervals that (simultaneously) enjoy the properties
exhibited in those Theorems.

\begin{theor}\label{laleche}
Let $L=(q_i)_{i\geq 1}$ be a sequence of non-negative, relatively prime integers. Let
$I=(a_j)_{j\geq 1}$ be a sequence such that every $a_j\in L$,
while each $q_i$ in $L$ appears infinitely many times in $I$. Let $J=(b_i)_{i\geq 1}$ be a strictly increasing sequence of non-negative integers such that $\gcd(q_{i},b_{j})=1$ for all
$i,j\geq 1$. Let $\gn (I)$ and $\gn (J)$ be the generalized
integers associated to $I$ and $J$ respectively. Then, for any
generalized integer $\gm$ coprime with $\gn (I)\cdot \gn (J)$,
there exists a simple Riesz group of rank one $G(\gm )$ such that:
\begin{enumerate}
\item For every $q_i$ in $L$, there is a countably generated
interval $D_i$ satisfying $tD_i\neq G(\gm )^+$ for $t\leq q_i-1$
and $q_i D_i=G(\gm )^+$. \item There is an interval $D\subset
G(\gm )^+$ such that $nD\neq G(\gm )^+$ for
 all $n$ in
$\N$.
\item For some generalized integer $\gn$ dividing $\gn (J)$,
the group $G(\gm)$ is isomorphic to $\Z_{\gn (I)\cdot \gn\cdot \gm}$ (as
abelian groups).
\end{enumerate}
\end{theor}
\begin{proof}
We first use the argument in the proof of Theorem~\ref{ertorema}
with the sequence $I$. In this way we get a simple group of rank
one $(H,H^+)$ such that: (a) $H\cong \Z_{\gn (I)}$; (b) For every
$q_i$ in $L$, there is a countably generated interval $E_i$
satisfying $tE_i\neq H^+$ for $t\leq q_i-1$ and $q_i E_i=H^+$. By
Proposition~\ref{Bergman 1}, $(H,H^+)=\varinjlim ((\Z, H_i^+),
l_i\cdot)$, where $(\Z, H_i^+)$ is a simple component and
$l_i\cdot\,\colon (\Z, H_i^+)\to (\Z, H_{i+1}^+)$ is an
order-embedding for all $i\geq 1$. Notice that $\prod_{i\geq
1}l_i=\gn (I)=\prod_{k\geq 1}q_k^{\infty}$. Thus, for each $i\geq
1$, $l_i=\prod_{j=1}^{r_i}{q_{k_j}^{n_j}}$ for some $r_i$ and
$n_j$ in $\N$. Therefore $\gcd (l_i, b_j)=1$ for all $i,j\geq 1$.

Now, fix $(\Z, K_1^+)=(\Z, H_1^+)$, and apply the argument in the
proof of Theorem~\ref{ertorema2}, using the sequence $J$. Thus, we
get an inductive system $((\Z, K_i^+), a_i\cdot)$, where $(\Z,
K_i^+)$ is a simple component, $a_i\in J$ and $a_i\cdot\,\colon
(\Z, K_i^+)\to (\Z, K_{i+1}^+)$ is an order-embedding for all
$i\geq 1$. Moreover, the group $(K,K^+)=\varinjlim ((\Z, K_i^+),
a_i\cdot)$ is a simple group of rank one such that: (a) $K\cong
\Z_{\gn}$ for the generalized integer $\gn =\prod_{n\geq 1}a_n$,
which divides $\gn (J)$; (b) There is a countably generated
interval $E$ such that $nE\ne K^+$ for every $n\geq 1$.

We next define submonoids $G_{i,j}^+$ of the non-negative integers
by recurrence on $i,j\geq 0$, as follows:
\begin{enumerate}
\item[(a)] $G_{0,0}^+=H_1^+=K_1^+$. \item[(b)] For every $i\geq
1$, $G_{i,0}^+=K_{i+1}^+$. \item[(c)] For every $j\geq 1$,
$G_{0,j}^+=H_{j+1}^+$. \item[(d)] For every $i,j\geq 1$,
$G_{i,j}^+=a_iG_{i-1,j}^++l_jG_{i,j-1}^+$.
\end{enumerate}

By~\cite[Lemma 2.3(1)]{P2}, we have that $(\Z, G_{i,j}^+)$ is a
simple component for every $i,j\geq 0$, and in the following
diagram:

\begin{equation}
\label{diagonales}
\begin{CD}
(\Z,G_{0,0}^{+}) @>{l_1\cdot}>> (\Z,G_{0,1}^{+})@>{l_2\cdot}>>
(\Z,G_{0,2}^{+})@>{l_3\cdot}>>\cdots\\ @V{a_1\cdot}VV
@VV{a_1\cdot}V @VV{a_1\cdot}V \\ (\Z,G_{1,0}^{+}) @>{l_1\cdot}>>
(\Z,G_{1,1}^{+})@>{l_2\cdot}>>
(\Z,G_{1,2}^{+})@>{l_3\cdot}>>\cdots\\ @V{a_2\cdot}VV
@VV{a_2\cdot}V @VV{a_2\cdot}V \\ (\Z,G_{2,0}^{+}) @>{l_1\cdot}>>
(\Z,G_{2,1}^{+})@>{l_2\cdot}>>
(\Z,G_{2,2}^{+})@>{l_3\cdot}>>\cdots\\ @V{a_3\cdot}VV
@VV{a_3\cdot}V @VV{a_3\cdot}V \\ \vdots  & & \vdots &  &
\vdots \\
\end{CD}
\end{equation}
all squares are commutative and all the maps are order-embeddings
(see Proposition~\ref{ladiagonal} and~\cite[Proposition
2.10]{P2}).

Let $(G,G^+)=\varinjlim ((\Z, G_{i,i}^+), a_il_i\cdot)$. Then
$(G,G^+)$ is a simple group of rank one, and $G\cong \Z_{\gn\cdot\gn
(I)}$ by construction. An argument analogous to that in condition
(iii) of Proposition~\ref{ladiagonal} guarantees the existence of
order-embeddings $\sigma\colon (H,H^+)\to (G,G^+)$ and $\tau\colon
(K,K^+)\to (G,G^+)$. Thus, for any generalized integer $\gm$
coprime with $\gn (I)\cdot\gn $, there exist by Theorem~\ref{2.11 nou}
a simple Riesz group of rank one $(G(\gm),G(\gm)^+)$ and an
order-embedding $\beta\colon (G,G^+)\to (G(\gm ), G(\gm )^+)$.
Clearly, the maps $(\beta\sigma)\colon (H,H^+)\to (G(\gm ), G(\gm
)^+)$ and $(\beta\tau)\colon (K,K^+)\to (G(\gm ), G(\gm )^+)$ are
order-embeddings. Thus, by condition (iii) in
Lemma~\ref{imageinterval}, the intervals $D=E_{(\beta\tau)}$ and
$D_i=(E_i)_{(\beta\sigma)}$ in $(G(\gm ), G(\gm )^+)$
enjoy the desired properties.
\end{proof}

The example in Theorem~\ref{laleche} above allows us to construct a (stably finite) monoid of intervals $W_{\sigma}^D(G^+)$ over a simple Riesz group $G$, where $D$ is an unbounded interval such that the representation map $\varphi$ defined in (\ref{lafi}) is not injective, even in the case when $D$ is not functionally complete (see~\cite[Remark 3.4(2)]{op}). Other consequences will be outlined in Section~7.

\section{Final comments and remarks}

In this section we explore the possible applications of the
results obtained in previous sections to the context of K-Theory
of multiplier algebras of simple $C^*$-algebras with real rank
zero.

We remind the reader that $C^*$-algebras are precisely the
norm-closed $^*$-subalgebras of $\mathbb{B}(\mathcal{H})$, the
algebra of bounded linear operators on a Hilbert space
$\mathcal{H}$. Recall that a (unital) $C^*$-algebra $A$ has
\emph{real rank zero} provided that the set of invertible
self-adjoint elements of $A$ is dense in the set of self-adjoint
elements of $A$ (see~\cite{B-P}). In case $A$ does not have a
unit, then $A$ has real rank zero if, by definition, the minimal
unitization $\widetilde{A}$ has real rank zero. We say that a
(unital) $C^*$-algebra $A$ has \emph{stable rank one} if the set
of invertible elements of $A$ is dense (see~\cite{Ri},~\cite{hv}).
As with the real rank zero case, if $A$ does not have a unit, then
$A$ has stable rank one if $\widetilde{A}$ has. A simple and
separable $C^*$-algebra is said to be \emph{elementary} if it is
isomorphic to the algebra of compact operators on a (separable)
Hilbert space. This translates into the requirement that the
algebra contains minimal projections. We shall be concerned with
non-elementary $C^*$-algebras.
\begin{probl}
\label{probl} \rm{Let $(G,G^+)$ be any of the groups
obtained in Theorems~\ref{ertorema},~\ref{ertorema2}
or~\ref{laleche}. Does there exist a simple, separable, non-unital
$C^*$-algebra $A$ with real rank zero and stable rank one for
which the ordered group $(K_0(A),K_0(A)^+)$ is order-isomorphic to
$(G,G^+)$ ?}
\end{probl}

We comment below on the relevance of this question, for the
consequences that would result given a positive answer. For this,
we need to remind the reader of some basic elements in K-Theory
that will be needed in our discussion (see, e.g.~\cite{Black}).
Given a $C^*$-algebra $A$, we denote by $M_{\infty }(A)=\varinjlim
M_n(A)$, under the maps $M_n(A)\to M_{n+1}(A)$ defined by
$x\mapsto \mbox{diag}(x,0)$, that is, $M_{\infty }(A)$ is the
algebra of countably infinite matrices over $A$ with only finitely
many non-zero entries.

We denote by $V(A)$ the set of equivalence classes of projections
in $M_{\infty}(A)$ under the Murray-von Neuman equivalence $\sim$.
This becomes an abelian monoid with operation $[p]+[q]=[\left(
\begin{array}{cc}
  p & 0 \\
  0 & q \\
\end{array}
\right)]$. This monoid is naturally endowed with the
\emph{algebraic pre-order}, denoted by $\leq$, induced by the
previous equivalence; namely $[p]\leq [q]$ if $p$ is equivalent to
a subprojection of $q$.

If the $C^*$-algebra $A$ is represented faithfully as a
$^*$-subalgebra of $\mathbb{B}(\mathcal{H})$, for a separable
Hilbert space $\mathcal{H}$, and the action is non-degenerate,
then we define the \emph{multiplier algebra} $\mathcal{M}(A)$ of
$A$ as the $C^*$-algebra
\[
\mathcal{M}(A)=\{x\in\mathbb{B}(\mathcal{H})\mid xA\subset A\mbox{
and }Ax\subset A\}\,.
\]
It is well-known that this construction is equivalent to the one
obtained by using double centralizers (see, e.g.~\cite{W-O}), and
it is of course only relevant in case $A$ does not have a unit
itself, since otherwise $\mathcal{M}(A)$ coincides with $A$. The
multiplier algebra, together with the embedding of $A$ as a
two-sided closed ideal, provides the solution to the universal
problem of adjoining a unit to the algebra $A$.

If $A$ is a separable (non-unital) $C^*$-algebra with real rank
zero and $P$ is a projection in $\mathcal{M}(A)$, then
by~\cite[Lemma 1.3]{K0Good} we have that $PAP$ also has real rank
zero and an approximate unit consisting of an increasing sequence
of projections, say $(p_n)$. If, moreover, $p$ is a projection in
$A$, then $p\lesssim P$ if and only if $p\lesssim p_n$ for some
$n\geq 1$. In this situation, we define
\begin{align*}
\Theta ([P]) & =\{[p]\in V(A)\mid p\mbox{ is a projection in } PM_{\infty }(A)P\}\\
 & =\{
[p]\in V(A)\mid [p]\leq [p_n] \mbox{ for some } n\mbox{ in
}\N\}\,.
\end{align*}
Then $\Theta ([P])$ is a countably generated interval in $V(A)$,
which is soft precisely when $P\notin A$. Let $D(A)=\Theta
([1_{\mathcal{M}(A)}])$. In the case that $A$ has moreover stable
rank one, then the map
\begin{equation}
\label{lateta} \Theta
:(V(\mathcal{M}(A)),[1_{\mathcal{M}(A)}])\rightarrow W_{\sigma
}^{D(A)}(V(A))
\end{equation}
is a normalized monoid isomorphism (see~\cite[Theorem 2.4]{Perera}
and also~\cite[Theorem 1.10]{K0Good}).

For any separable, non-unital, non-elementary simple $C^*$-algebra
with real rank zero and stable rank one, it is well-known that the
group $K_0(A)$ is a countable, non-atomic, simple, Riesz group.
Because of the existence of an approximate unit of projections,
$K_0(A)$ is naturally isomorphic to the Grothendieck group of the
monoid $V(A)$. Since $A$ has stable rank one, $V(A)$ has
cancellation and it can be identified with $K_0(A)^+$. Let $p$ be
any non-zero projection in $A$ and set $u=[p]$ in $V(A)$. If
$d=\sup\phi_u (D(A))$ (see also the notation in Section~3), then
by composing the map $\varphi$ defined in~(\ref{lafi}) with the
map defined in~(\ref{lateta}), we get a normalized monoid morphism
\begin{equation}
\label{laPhi} \Phi\colon (V(\mathcal{M}(A)),
[1_{\mathcal{M}(A)}])\to (V(A)\sqcup W_{\sigma }^d(S_{u}), d)\,,
\end{equation}
which is an isomorphism if $V(A)$ is furthermore strictly
unperforated, see~\cite[Theorem 3.8]{Perera}. We now comment on
what kind of examples a positive solution to Problem~\ref{probl}
would lead to in connection with the results obtained in previous
sections.
\begin{p}
\label{toma1} Let $I=(q_n)$ be an increasing sequence of
non-negative and relatively prime integers, and let $A$ be a
separable, non-unital, non-elementary $C^*$-algebra with real rank
zero and stable rank one such that $(K_0(A), K_0(A)^+)$ is
order-isomorphic to the group constructed in
Theorem~\ref{ertorema} (with respect to the sequence $I$). Then,
there exists a sequence of projections $(E_n)_{n\geq 1}\subset
\mathcal{M}(A\otimes\mathbb{K})$ such that:
\begin{enumerate}
\item $\Phi (E_n)=\Phi
(1_{\mathcal{M}(A\otimes\mathbb{K})})=\infty$ for all $n\geq 1$.
\item $E_{n+1}\lnsim E_n$ for every $n\geq 1$ (i.e. $E_{n+1}$ is not equivalent to a subprojection of $E_n$). \item For each
$n\geq 1$, $t\cdot E_n\nsim 1_{\mathcal{M}(A\otimes\mathbb{K})}$
whenever $t< q_n$, and $q_n\cdot E_n\sim
1_{\mathcal{M}(A\otimes\mathbb{K})}$.
\end{enumerate}
\end{p}

Replace $A$ by its stabilization $A\otimes\mathbb{K}$ (where
$\mathbb{K}$ is the $C^*$-algebra of compact operators on a
separable Hilbert space), and note that the $K_0$ group remains
the same. So, to verify the above claim, assume that $A$ is stable.

By Theorem~\ref{ertorema}, for every $q_{n}$ in $I$, there is a
countably generated unbounded interval $D_n\subset K_0^+(A)$ such
that $tD_n\neq K_0^+(A)$ for $t\leq q_{n}-1$ and $q_n
D_n=K_0^+(A)$. Moreover, by Theorem~\ref{eldarrer}, we can choose
these intervals in such way that $D_{n+1}\leq D_n$ in the
algebraic ordering of the monoid of intervals
$W_{\sigma}^{D(A)}(K_0^+(A))$.

Since, as mentioned, we can identify $V(A)$ with $K_0(A)^+$, we
can use the isomorphism~(\ref{lateta}), to get a sequence of
projections in $\mathcal{M}(A\otimes\mathbb{K})$ by setting
$E_n=\Theta^{-1} (D_n)$. Clearly they satisfy properties
(i)-(iii).

Notice that, if $A$ is a $C^*$-algebra satisfying the hypotheses
in~\ref{toma1} then, for every $n$, the $C^*$-algebra
$\mathcal{M}(E_n(A\otimes \mathbb{K})E_n)$ is finite. Otherwise,
at the level of monoids, $D_n+Y=D_n$ for a non-zero interval $Y$,
and thus, by simplicity of $K_0(A)$, we would conclude that
$D_n=K_0^+(A)$, in contradiction with Theorem~\ref{ertorema}. On
the other hand
\[
M_{q_n}(\mathcal{M}(E_n(A\otimes \mathbb{K})E_n))\cong
\mathcal{M}(M_{q_n}(E_n(A\otimes \mathbb{K})E_n))\cong
\mathcal{M}(A\otimes \mathbb{K})\,,\] which implies that
$M_{q_n}(\mathcal{M}(E_n(A\otimes \mathbb{K})E_n))$ is not finite.
This kind of behaviour has been exhibited in concrete examples
constructed by R\o rdam (see~\cite{Rordam}). There are even simple
examples, but they don't have real rank zero
(see~\cite{ror2},~\cite{ror3}).

The existence of examples as in~\ref{toma1} would provide us with examples of $C^*$-algebras of real rank zero that fail to have weak cancellation in the sense of Brown and Pedersen (see~\cite{bppreprint}). They would also give a solution to the Fundamental Separativity Problem (see, e.g.~\cite{agop}).

\begin{p}
\label{toma2} Let $A$ be a separable, non-unital, non-elementary
$C^*$-algebra with real rank zero and stable rank one such that
$(K_0(A), K_0(A)^+)$ is order-isomorphic to the group constructed
in Theorem~\ref{ertorema2}. There exists then a projection $E$ in
$\mathcal{M}(A\otimes\mathbb{K})$ such that:
\begin{enumerate}
\item $\Phi (E)=\Phi
(1_{\mathcal{M}(A\otimes\mathbb{K})})=\infty$. \item $n\cdot
E\nsim 1_{\mathcal{M}(A\otimes\mathbb{K})}$ for every $n\geq 1$.
\end{enumerate}
\end{p}

To check this, use Theorem~\ref{ertorema2}, to find a countably
generated unbounded interval $D\subset K_0^+(A)$ such that $nD\neq
K_0^+(A)$ for every $n\geq 1$. Then, using the
isomorphism~(\ref{lateta}), we get a projection $E=\Theta^{-1}
(D)$ in $\mathcal{M}(A\otimes\mathbb{K})$ satisfying the required
properties.

Notice that if $A$ is a $C^*$-algebra satisfying the hypotheses
in~\ref{toma2}, then we have an answer to an implicit question
posed in~\cite[Remark 3.4(2)]{op}. Namely, if $(G,G^+)$ is a
simple Riesz group containing an interval $D\subseteq G^+$ such
that $\varphi (D)=\varphi (G^+)=\infty$, but $nD\ne G^+$ for every
$n$ in $\N$, then if $Y\in W_{\sigma
}^{D}(G^+)$ and $nD+Y=nD$, the simplicity of $(G, G^+)$ would imply
that $(n+1)D=nD$, and thus $nD=G^+$, contradicting the
hypothesis (i.e. $W_{\sigma }^{D}(G^+)$ is a stably finite monoid,
see e.g.~\cite{op}). So, $nD\ne mD$ whenever $n\ne m$, but still
$\varphi (nD)= \infty$. Hence, it might be possible to construct a
unital, simple $C^*$-algebra $A$ with real rank zero and stable
rank one, such that the multiplier algebra $\mathcal{M}(A\otimes
\mathbb{K})$ contains a non-zero projection $E$ with
$\mathcal{M}(E(A\otimes \mathbb{K})E)$ stably finite, but with
identically infinite scale~(\cite{Perera}). Moreover, according
to~\cite[Proposition 3.6]{Rordam} (also see~\cite[Theorem
2.10]{Pardmult}), $E(A\otimes \mathbb{K})E$ would not be a stable
algebra. The existence of such an example would fix the exact
limits of application of~\cite[Proposition 2.11]{Pardmult}.

\begin{p}
\label{toma3} Let $I=(q_n)$ be an increasing sequence of
relatively prime non-negative integers. Let $A$ be a separable,
non-unital, non-elementary $C^*$-algebra with real rank zero and
stable rank one such that $(K_0(A), K_0(A)^+)$ is order-isomorphic
to the group constructed in Theorem~\ref{laleche} (with respect to
the sequence $I$). There exists then a sequence of projections
$(E_n)_{n\geq 1}$ and a projection $E$ in
$\mathcal{M}(A\otimes\mathbb{K})$ such that:
\begin{enumerate}
\item $\Phi (E_n)=\Phi
(1_{\mathcal{M}(A\otimes\mathbb{K})})=\infty$ for every $n\geq 1$.
\item $E_{n+1}\lnsim E_n$ for every $n\geq 1$. \item For every
$n\geq 1$, $t\cdot E_n\nsim 1_{\mathcal{M}(A\otimes\mathbb{K})}$
whenever $t< q_n$, and $q_n\cdot E_n\sim
1_{\mathcal{M}(A\otimes\mathbb{K})}$. \item $\Phi (E)=\Phi
(1_{\mathcal{M}(A\otimes\mathbb{K})})=\infty$. \item $n\cdot
E\nsim 1_{\mathcal{M}(A\otimes\mathbb{K})}$ for every $n\geq 1$.
\end{enumerate}
\end{p}

Hence, in view of~\ref{toma1} and~\ref{toma2}, the existence of a
$C^*$-algebra $A$ satisfying the hypotheses in~\ref{toma3} would
imply that the multiplier algebra $\mathcal{M}(A\otimes
\mathbb{K})$ contains projections having the special behaviors
stated in there.

\section*{Acknowledgments}

Part of this work was done during a visit of the first author to the
Department of Pure Mathematics at Queen's University Belfast
(Northern Ireland), of the second author to the Centre de Recerca
Matem\`atica, Institut d'Estudis Catalans in Barcelona (Spain), and
of the third author to the Departamento de Matem\'aticas de la
Universidad de C\'adiz (Spain). The three authors are very grateful
to the host centers for their warm hospitality.

\end{document}